\documentclass[12pt,reqno]{amsart}
\usepackage{amsthm}
\usepackage{fullpage,url}
\usepackage[dvips]{hyperref}
\usepackage{amscd}   
\usepackage[all,graph,poly]{xy} 
\usepackage{amssymb,amsmath}
\usepackage{graphicx}
\usepackage{pstricks,pst-plot,pst-node}
\usepackage[dvips]{hyperref}
\usepackage{multicol}
\usepackage{setspace}
\DeclareFontEncoding{OT2}{}{} 

\providecommand{\abs}[1]{\lvert#1\rvert}

\newcommand{\cF}{\mathcal{F}}

\newcommand{\op}[1]{\operatorname{#1}}

\newcommand{\msf}[1]{\mathsf{#1}}

\newcommand{\ZZ}{\mathbb{Z}}

\DeclareMathOperator{\cl}{cl}

\DeclareMathOperator{\rk}{rk}
\DeclareMathOperator{\st}{star}

\def\co{\colon\thinspace}
\newcommand{\e}[1]{\ar@{-}[#1]}
\newcommand{\ed}[1]{\ar@{--}[#1]}
\newcommand{\ee}[1]{\ar@{=}[#1]}
\newcommand{\eer}[1]{\ar@{=}[#1] |{\SelectTips{eu}{}\object@{<}} |>{\SelectTips{}{}\object@{}}}
\newcommand{\eel}[1]{\ar@{=}[#1] |{\SelectTips{eu}{}\object@{>}} |>{\SelectTips{}{}\object@{}}}
\newcommand{\er}[1]{\ar@{-}[#1] |{\SelectTips{eu}{}\object@{<}}| {\SelectTips{eu}{}\object@{}}}
\newcommand{\el}[1]{\ar@{-}[#1] |{\SelectTips{eu}{}\object@{>}}{\SelectTips{eu}{}\object@{}}}
\newcommand{\edr}[1]{\ar@{--}[#1] |{\SelectTips{eu}{}\object@{<}} |{\SelectTips{eu}{}\object@{}}}
\newcommand{\edl}[1]{\ar@{--}[#1] |{\SelectTips{eu}{}\object@{>}} |{\SelectTips{eu}{}\object@{}}}
\newcommand{\ert}[1]{\ar@{-}[#1] |---------{\SelectTips{cm}{}\object@{>}}|>{\SelectTips{eu}{}\object@{}}}
\newcommand{\eroff}[1]{\ar@{-}[#1] |--------{\SelectTips{eu}{}\object@{<}}|>{\SelectTips{eu}{}\object@{}}}
\newcommand{\eloff}[1]{\ar@{-}[#1] |--------{\SelectTips{eu}{}\object@{>}}|>{\SelectTips{eu}{}\object@{}}}
\newcommand{\eeroff}[1]{\ar@{=}[#1] |---------{\SelectTips{eu}{}\object@{<}} |>{\SelectTips{}{}\object@{}}}
\newcommand{\eeloff}[1]{\ar@{=}[#1] |---------{\SelectTips{eu}{}\object@{>}} |>{\SelectTips{}{}\object@{}}}

 \newcommand{\lul}[1]{\ar@{}[l]_<<{#1}}
\newcommand{\rrul}[1]{\ar@{}[r]^<<<<{#1}}
\newcommand{\rul}[1]{\ar@{}[r]^<<{#1}}
\newcommand{\ldl}[1]{\ar@{}[l]^<<{#1}}
\newcommand{\rdl}[1]{\ar@{}[r]_<<{#1}}
\newcommand{\dl}[1]{\ar@{}[d]_<<{#1}}
\newcommand{\dll}[1]{\ar@{}[dd]_{#1}}
\newcommand{\rrdl}[1]{\ar@{}[ru]_<<{#1}}
\swapnumbers
\newtheorem{theorem}{Theorem}[section]
\newtheorem{lemma}[theorem]{Lemma}

\newtheorem{proposition}[theorem]{Proposition}
\newtheorem{lemma/definition}[theorem]{Lemma/Definition}

\theoremstyle{definition}
\newtheorem{definition}[theorem]{Definition}
\newtheorem{definition/lemma}[theorem]{Definition/Lemma}
\newtheorem{definition/remark}[theorem]{Definition/Remark}

\newtheorem{example}[theorem]{Example}

\newtheorem{topic}[theorem]{}
\theoremstyle{remark}
\newtheorem{remark}[theorem]{Remark}
\begin{document}
%
%
%
%
\title{Combinatorial cell complexes and Poincare duality}
\author[Tathagata Basak]{Tathagata Basak} 
\address{Department of Mathematics\\ University of Chicago\\ Chicago, IL 60637\\USA}
\email{tathagat@math.uchicago.edu}
\urladdr{http://www.math.uchicago.edu/~tathagat}
\date{July 25, 2008}
\begin{abstract}
We define and study a class of finite topological spaces, which model the cell
structure of a space obtained by gluing finitely many Euclidean convex polyhedral
cells along congruent faces. 
We call these finite topological spaces, combinatorial cell complexes (or c.c.c).
We define orientability, homology and cohomology of c.c.c's and develop enough
algebraic topology in this setting to prove the Poincare duality theorem for
a c.c.c satisfying suitable regularity conditions. The definitions and proofs
are completely finitary and combinatorial in nature.
\end{abstract}
\maketitle
%
%
%
%
\section{Introduction}
%
%
%
%
\begin{topic}{\bf Summary of results: }
Given a topological space with a triangulation, if we only remember the set of
simplices and incidence relations among them, we get a simplicial complex.
One can think of the partially ordered set of the simplicial complex as a finite
topological space and study how the combinatorics of this poset reflects the
algebraic topology of the space one started with.
In this article we want to do something similar, but we want to allow
our cells to have more general shapes, not just of simplices. (For example, cells in
the shape of any convex polyhedron are allowed).
We shall call these objects combinatorial cell complex or c.c.c for short.
Let $X$ be a topological space written as a finite union of a collection
$S_X$ of Euclidean convex polyhedra.
Assume that $S_X$ is closed under intersection and that the intersection of two distinct
polyhedron in $S_X$ of equal dimension, has strictly lower dimension.
If we forget the  space $X$ and only remember the set $S_X$, the dimension of each
polyhedron and the partial order coming from incidence relation among the elements
of $S_X$, we get an example of a c.c.c.
\par
 Thus, a c.c.c $S$ is a partially ordered set, with a rank (or dimension) function defined
on $S$, satisfying some axioms (the definition is given in \ref{cccdef}).
The elements of $S$ are called cells. The axioms describe how the cells
are allowed to be glued together; they try to mimic the  conditions that are satisfied if
$S$ was obtained from polyhedral decomposition of a space $X$, as above.
Our objective here are the following:
\newline
(A) We want to see how to translate into $S_X$,
the topological properties of $X$, via the correspondence $X \to S_X$.
For example, we shall call $S$ manifold--like, if it satisfies some extra conditions
that would obviously hold, if $S = S_X$ for some manifold $X$.
The key definition is that of an orientable c.c.c (see \ref{orientationdef}).
\newline
(B) Once we have put enough regularity conditions on a c.c.c to remove the pathologies,
we want to see how much algebraic topology can be developed
in this combinatorial setting. In particular, we define cellular homology and cohomology
groups of c.c.c's with orientable cells and prove a Poincare duality theorem stated below
(see theorem \ref{duality}).
\par
{\bf Theorem. }{\it Let $S$ be an orientable,  manifold--like c.c.c of dimension $n$.
Suppose each cell of the c.c.c $S$ and the opposite c.c.c $S^{\circ}$ is
flag--connected and acyclic. Then $H_i( S) \simeq H^{n-i}(S)$.}
\par
(The definitions of the various terms are given in the following sections:
flag--connected and orientable: \ref{orientationdef}, manifold--like: \ref{manifolddef},
$S^{\circ}$: \ref{Sop}, acyclic: \ref{acyclicdef}.
Homology and cohomology groups  are defined in section \ref{sechomcohom}).
If $S = S_X$ for some space $X$, then these homology groups are the same as
the cellular homology groups of $X$. In particular, a simplicial complex
gives a c.c.c and in this situation, our homology groups are identical with
simplicial homology groups (see \ref{simphom}). 
\par
The main technical part in the proof of theorem \ref{duality} is to show that,
under the conditions of the theorem, the homology of $S$ 
is invariant under ``barycentric subdivision" (see proposition \ref{baryinv}).
It follows (see \ref{functor}), that under these regularity conditions, the homology groups
of the c.c.c $S$ coincide with the homology groups of the simplicial set $N(S)$ obtained by
taking the nerve of the poset $S$ (or, in other words, the singular homology of the topological
space obtained by taking geometric realization of $N(S)$). Sections \ref{secsteller}, 
\ref{secvanishing} and \ref{secbarycentric} are mainly occupied with proving \ref{baryinv}.
Given the technical result \ref{baryinv}, the proof of the theorem \ref{duality} is totally
transparent. This argument, given in section \ref{secthm}, can be read right after we are
through with the definitions in section \ref{sechomcohom}. 
\end{topic}
%
%
\begin{topic}{\bf Relationship with simplicial topology: }
The standard approach for translating algebraic topology in a combinatorial setting
is via simplicial sets (e.g. see \cite{PM:SA}), which are abstract versions of simplicial complexes
with labeling of vertices. Our main reason for introducing a combinatorial setting with more
general cell shapes is the following:
\par
In the classical proof of Poincare duality, one relates homology and cohomology by
taking the dual of a cell complex (e.g. see \cite{GH:AG}).
However, the cells of the dual cell complex of a simplicial complex need not be simplices.
We allow more general cell shapes so that the duality is built into
the setup (the dual of a c.c.c $S$ has the same underlying set as $S$, with the partial
order and rank reversed).
\par
One disadvantage of the present setup is the lack of explicit functoriality of homology
groups. In general, given a continuous map (that is, an order preserving function) $f \co S \to S'$
between c.c.c's, there is no obvious chain map from the chain complex of $S$
(as defined in section \ref{sechomcohom}) to that of $S'$, inducing a map between the cellular
homology groups.
However, if $S$ and $S'$ satisfy the regularity conditions given in the Poincare duality
theorem above, then one does get a map $H_i(f) \co H_i(S) \to H_i(S')$, so that
$H_i$ becomes a functor. Unfortunately, we are only able to see this by using the 
invariance of homology under barycentric subdivision (see \ref{baryinv}, \ref{baryinvf}), 
 and the consequent canonical isomorphism between the homology of a
c.c.c $S$ (with enough regularity conditions) and that of the simplicial set $N(S)$
(see \ref{functor}).
The functoriality of the cellular homology groups follows by invoking
the functoriality of homology of simplicial sets.
\par
As was suggested by Peter May (private communications), it would be nice to have
a shape category so that (some variant of) a c.c.c becomes a presheaf (of sets) on
this shape category. Then one could develop the theory as for simplicial sets in a
functorial way.  This possibility also makes us wonder if the combinatorial study
of shapes of cells might have some bearings on certain approaches to higher category theory,
notably those initiated by Steet in \cite{RS:AS} and by Baez--Dolan in \cite{BD:HA}.
In these approaches, much of the structure of the higher category is encoded in
the shape of the cells that represent the higher morphisms. For an introduction to these
ideas, see chapter 6 and 4 in \cite{CL:HC}.
\end{topic}
%
%
\begin{topic}{\bf Finite topological spaces: }
The topology of finite spaces can be surprisingly rich. For example, there are finite spaces
having weak homotopy type of any finite simplicial complex (see \cite{MM:SH}). The
finite topological spaces have been studied since they were introduced by Alexandroff
in \cite{PA:DR} and the theory of simple homotopy types was developed by Whitehead in \cite{JW:SH}.
The simple homotopy types of polyhedra were studied using finite topological spaces
in the recent article \cite{BM:SF}. 
We refer the reader to the notes \cite{PM:FT} and \cite{PM:FS} for an introduction to the
topology of finite spaces and to \cite{MW:PT} for a survey of the combinatorial aspects of this
theory. The book \cite{DK:CA} is a convenient reference for combinatorial algebraic topology.
\par
In this article we have restricted our study to purely combinatorial aspects of the theory of c.c.c's.
The close relationship between the topology of a cell complex and that of the corresponding
finite space has not been explored or utilized here. This, and other topological questions,
like the relationship between the homology of a c.c.c $S$ defined here and the singular
homology of the finite space $S$, will hopefully be explored in a later article.
\end{topic}
%
%
\begin{topic}{\bf Organization of the paper: }
The arguments in this article are, in most places, logically self contained.
The proof of some technical lemmas have been relegated to an appendix
to arrive at the main theorem \ref{duality} quickly.
An index of some frequently used symbols is included below.
\newline
\newline
{\bf Acknowledgments:}  I would like to thank Prof. Gabriel Minian for many useful comments on reading
an early draft of this article.  I would like to thank Prof. Richard Borcherds and Prof. Jon
Alperin for suggesting useful references. Most of all, I would like to thank Prof. Peter
May for his encouragement and many interesting and illuminating discussions since the early
stages of this work.
\end{topic}
%
%
\begin{topic}
{\bf Index of some commonly used notation: }
Let $S$ be a c.c.c. Let $T$ be a subset of $S$ and $x$, $y$ be elements of $S$.
\begin{tabbing}
$Delta x$ X\= the set of faces of $x$.X\=\kill
$C_i(S)$ \> the set of $i$--chains in the c.c.c $S$, that is, the free abelian group on the $i$--cells of $S$.\\
$C_x(y)$ \>  a ``new cell" in the stellar subdivision $S^x$, called the cone on $y$ with vertex at $x$.\\
$\op{cl}(x)$ \> the set of cells less than or equal to $x$, that is, the closure of $x$.\\
$\Delta(x)$ \> the set of faces of $x$.\\
$\partial$ \> the boundary map on chain complexes.\\
$\cF(S)$ \> the set of flags in $S$. (We write $\cF(x) = \cF(\cl(x))$). \\
$\gamma$ \> usually a flag (except in lemma \ref{simphom}, where it is a simplex).\\
$M(x)$      \> $= \st(x) \setminus U(x)$. \\
$\nabla(x)$ \> the set of co-faces of $x$.\\
$\omega$ \> an orientation. ($\omega_x$ denotes an orientation on $\cl(x)$).\\
$\rk(x)$ \> the rank of a cell $x$.\\
$S$ \> usually a combinatorial cell complex (called c.c.c for short). \\
$S^{\circ}$ \> the opposite c.c.c of $S$.\\
$S(r)$ \> the set of cells of $S$ having rank $r$. \\
$S^{(1)}$ \> the first barycentric subdivision of $S$. \\
$S^x$ \> the stellar refinement of $S$ at $x$. We shall write $(S^x)^y = S^{x y}$.\\
$s(x,y)$ \> a sign assigned to each pair $\lbrace x, y \rbrace$,
where $x$ is a cell and $y$ is a face of $x$, with\\
         \> orientations given on $\cl(x)$ and $\cl(y)$ (see \ref{defsigns}).\\
$\st(x)$ \> $= \op{cl}(U(x))$, that is, the set of cells $y$ such that $x$ and $y$ have an upper bound.\\
$U(x)$ \> the set of cells greater than or equal to $x$.\\
$\vee T$ \> the least upper bound of $T$; we write $x \vee y = \vee \lbrace x, y \rbrace$. \\
$\wedge T$ \> the greatest lower bound of $T$; we write $x \wedge y = \wedge \lbrace x, y \rbrace$.\\
$X$ \> usually a combinatorial cell complex  (called c.c.c for short).\\
$x,y,z$ \> usually any of these letters denote a cell of a c.c.c.\\
\end{tabbing}
\end{topic}
%
%
%
%
\section{basic definitions}
%
%
%
%
\begin{topic}{\bf The setup: }
Suppose we are given a finite partially ordered set $(S, \leq)$
and a function, denoted by $\rk$, from $S$ to non-negative integers such that $y < x$ implies
$\rk(y) < \rk(x)$. Given this data, we introduce the following notation and nomenclature:
\par
If there is a possibility of confusion, we shall write $\leq_S$ to denote the partial order on $S$.
Elements of $S$ are called {\it cells}. If $\rk(x)= r$, we say that $x$ is a cell of
{\it rank} $r$ or $x$ is an $r$--cell. Write $S$ as a disjoint union,
$S = \cup_{r=0}^{\infty} S(r)$, where $S(r)$ is the set of $r$--cells of $S$.
If $x > y,$ we say that $x$ is above $y$ or that $y$ is below $x$. More precisely, we say
that $y$ is a {\it facet} of $x$ of co-dimension $(\rk(x) - \rk(y))$.
A co-dimension one facet of $x$ is
called a {\it face} of $x$. The set of faces of $x$ is denoted by $\Delta_S x $ or $\Delta x$, if
there is no possibility of confusion.
Dually, the cells that have $x$ as one of their face are called the {\it co-faces} of $x$.
The set of co-faces of $x$ is denoted by $\nabla x$.
The set of cells greater than or equal to $x$ (resp. less than or equal to $x$) is denoted by $U_S(x) = U(x)$
(resp. $\cl_S(x) = \cl(x)$).
\par
Let $T$ be a subset of $S$.
An element $x \in S$ is an {\it upper bound} of $T$, if $x \geq z$ for all $z \in T$.
The {\it least upper bound} of $T$, denoted by $\vee T$, is an upper bound 
of $T$ such that $\vee T \leq y$, for every upper bound $y$ of $T$.
Similarly, one defines the {\it greatest lower bound} of $T$, denoted by $\wedge T$.
Of course least upper bound or greatest lower bound of $T$ may not exist.
One also writes $x \vee y$ to denote $\vee \lbrace x, y \rbrace$ and $x \wedge y$ to denote
$\wedge \lbrace x, y \rbrace$. If $z = x \wedge y$, we say that $x$ and $y$ {\it meet} at $z$.
For $T \subseteq S$, let $\Delta T = \cup_{x \in T} \Delta x$.
Inductively define $\Delta^j T = \Delta ( \Delta^{j-1} T )$.
The rank zero cells below $x$ are called the {\it vertices} of $x$.
\end{topic}
%
%
\begin{definition}
\label{cccdef}
We say that $S$ is a {\it combinatorial cell complex} or c.c.c for short, if the data
$(S, \leq, \rk)$ satisfies the following four axioms:
\begin{enumerate}
\item  The partial order is compatible with rank, that is, if $y < x$,
then $\rk(y) < \rk(x)$.
\item The collection $S$ has enough cells, in the following sense.
If $T$ is a subset of $S$ that is bounded below, then the greatest lower bound $\wedge T$ exists.
For all $x$ and $y$ in $S$ with $y < x $, there exists a cell $y'$ such that $\rk(y') = \rk(y) + 1$ and
$y < y' \leq x$.
\item Each cell $x \in S$ of rank at-least one is the least upper bound of its faces, that is, $x = \vee \Delta x $.
\item  If $y$ is a co-dimension $2$ facet of $x$, then
there are exactly two faces of $x$ that are above $y$ and these two cells meet at $y$.
In other words, given $y \in S(i-1)$, $x \in S(i+1)$, $y < x$, there exists
distinct cells $y_+ $ and  $y_-$ in $S(i)$
such that $\Delta x \cap \nabla  y  =  \lbrace y_{+} , y_{-} \rbrace $.
\end{enumerate}
\end{definition}
%
%
\begin{example}
Let $T$ be an finite abstract simplicial complex (see definition 2.1 in \cite{DK:CA}).
The set $T$ becomes a combinatorial cell complex
with the partial order given by set inclusion. A simplex with $(r+1)$ vertices has
rank $r$. Given a collection of simplices $T_1 \subseteq T$, that is bounded below, the
greatest lower bound of $T_1$ is $\wedge T_1 = \cap_{\sigma \in T} \sigma$.
A co-dimension 2 facet of a simplex $\sigma$ has the form $\sigma \setminus \lbrace x_i, x_j\rbrace$,
where $x_i \neq x_j$ are two vertices of $\sigma$. The two simplices in between, are
$\sigma \setminus \lbrace x_i\rbrace$ and $\sigma \setminus \lbrace x_j \rbrace$.
\par
A topological space with a polyhedral decomposition defines a combinatorial cell complex.
Note that an $r$--cell has at-least $(r+1)$ vertices, but it can have more vertices.
\par
One can construct new combinatorial cell complexes from old ones by taking
sub-complexes (see \ref{subcomplex}), finite products
\footnote{Let $X$ and $Y$ be c.c.c's. The Cartesian product $X \times Y$ is a c.c.c,
with the induced partial order
(that is, $(x,y) \leq (x',y')$ if and only if $x \leq x'$ and $y \leq y'$) and rank given by
$\rk(x,y) = \rk(x) + \rk(y)$.},
barycentric and stellar subdivisions (see \ref{orientationdef} and \ref{stellerdef} respectively). 
\end{example}
%
%
\begin{topic}{\bf Topology on a c.c.c: }
Declare a subset $C$ of $S$ to be closed if $x \in C$ and $ y \leq x $ implies
$y \in C$. This defines a topology on $S$ in which arbitrary union and
intersection of closed sets are closed. Such spaces are called an A-space in \cite{PM:FT}.
({\it Caution: } What we are calling an closed set here is called an open set in
\cite{PM:FT} and vice versa. Both these conventions are found in the literature.)
Let $T$ be a subset of a c.c.c $S$. The closure of $T$, denoted by 
$\cl_S(T) = \cl(T)$, is the set of cells that are less than or equal to
some cell in $T$; these are precisely the closed sets of $S$.
If $x \in S$, then $\cl(x) = \cl(\lbrace x \rbrace)$ is the smallest closed set
containing $x$, so each cell of rank atleast one is a non-closed point in the
above topology. So $S$ is almost never Hausdorff. However $S$ is a $T_0$ space.
The subset $\lbrace x \in S \colon \rk(x) \leq i \rbrace$ is a closed subset of
$S$, called the {\it $i$-skeleton} of $S$.
\end{topic}
%
%
\begin{lemma}
(a) Let $C$ be a closed subset of $S$. Then $C$, with the rank and partial order
induced from $S$, is a c.c.c. \\
(b) Let $T \subseteq S$. Then the set of lower bounds of $T$ is equal to $\cap_{ x \in T} \cl (x) = \cl( \wedge T)$,
with the convention that $\cl(\wedge T) = \emptyset$, if $T$ is not bounded below.
\label{subcomplex}
\end{lemma}
\begin{proof}
(a) Axiom (1) holds for $C$ since the rank and partial order on $C$ are induced from $S$.
For axioms (2) and (4), we just need to observe that if $x \in C$ and $y \leq x$,
then $y \in C$. It also follows from this observation that $\Delta_C  x  = \Delta_S  x $, for all $ x \in C$.
This implies axiom (3), that is, $ \vee \Delta_C x  = x$. Part (b) follows from the definitions.
\end{proof}
\begin{remark} We end this section with a couple of easy observations.
The first one will be often used without explicit reference.
\begin{enumerate}
\item If $ z_+ \neq z_-$ are two cells with a common face $z$, then $z_+ \wedge z_- = z$.
So, if $x$ is a cell such that $z_+ > x$ and $z_- > x$, then $z = z_+ \wedge z_- \geq x$. Stated
differently, if $z \notin U(x)$, then at-most one of the co-faces of $z$ can belong to $U(x)$.
\item A subset $U$ of $S$ is open if and only if $x \in U$ and $ y \geq x$ implies $ y \in U$.
Thus $U(x) = \lbrace  y \in S \colon y \geq x \rbrace$ is the smallest open set containing $x$.
Given posets $S$ and $S'$, a function $f \co S \to S'$ is continuous in the above topology
if and only if it preserves the partial order. 
\end{enumerate}
\end{remark}
%
%
%
%
\section{nonsingular and manifold--like c.c.c.}
%
%
%
%
\begin{definition/remark}
\label{manifolddef}
A cell of a c.c.c is  maximal, if it is not below any other cell.
The {\it dimension} of a c.c.c $S$ is defined to be the maximal rank of a cell in $S$.
We say that $S$ is {\it equidimensional}, of dimension $n$, if each maximal
cell in $S$ has rank $n$.
\par
Assume that $S$ is {\it equidimensional}, of dimension $n$. The {\it boundary} of $S$
is defined  to be the set of cells of rank strictly less than
$n$, that have only one maximal cell above them. Since every cell of rank atleast one,
is the least upper bound of its faces, an $1$--cell cannot have only one vertex. So the {\it co-boundary}
of $S$, that is $ \lbrace y \in S(1) \colon \abs{ \Delta y } = 1 \rbrace $, is empty.
\par
A c.c.c $S$ of dimension $n$ is called {\it non-singular} if $S$ is equidimensional, each
$(n-1)$--cell of $S$ is a face of at-most two maximal cells and dually, each
$1$--cell of $S$ has at-most two vertices (hence exactly two vertices).
\par
We say that $S$ is {\it manifold--like} if it is nonsingular and has empty boundary.
Axiom (4) in definition \ref{cccdef} implies that the boundary of $\cl( \Delta x) $ is empty
for all $x \in S$.
\end{definition/remark}
%
%
\begin{lemma} Let $S$ be a c.c.c.\\
(a) For each $x \in S(r)$ and $ 0 \leq j \leq r$ one has,
\newline
(i) $\Delta^j  x  = \lbrace y \in S(r -j) : y \leq x \rbrace$,
\newline
(ii) $\vee \Delta^j  x  = x $. 
\newline
(b) Every subset of $S$, that is bounded above, has a least upper bound.\\
(c) Let $S$ be manifold--like, of dimension $n$ and $x \in S(r)$ for some $r < n$.
Then $\wedge \nabla  x  = x$.\\
(d) For all $x < y$ in $S$, one has $\Delta y \nsubseteq U(x)$.
\label{firstlemma}
\end{lemma}
\begin{proof}
(a) Axiom (2) in definition \ref{cccdef} implies that a co-dimension $j$ facet of
$x$ is a face of a co-dimension $(j-1)$ facet. The statement (i) follows from this by
induction on $j$.
\par
The proof of (ii) is also by induction on $j$. The case $ j = 1$ is the axiom (3) in definition
\ref{cccdef}. Notice that axiom (2) in definition \ref{cccdef} has the following consequence: if
$z_j \in \Delta^j  x $, then there exists $z_j < z_{j-1} < \dotsb < z_1 < z_0 = x$
such that $z_r$ is a facet of $x$ of co-dimension $r$.
It follows that $\Delta^j x  = \cup_{ y \in \Delta x} \Delta^{j-1} y$.
By induction, we may assume that $ \vee  \Delta^{j-1} y = y$.
Clearly $x$ is an upper bound for $\Delta^j x$.
Let $u$ be any upper bound of $\Delta^j x$. Then $u \geq t$ for all $t \in \Delta^{j-1} y$
and for all $y \in \Delta x$. Hence $u \geq \vee \Delta^{j-1} y = y$ for each $y \in \Delta x$.
It follows that $u \geq  \vee \lbrace y \colon y \in \Delta x \rbrace = x$.
\par
(b) If the set of upper bounds of $T$ is non-empty, it is easy to see that the greatest lower bound of
the upper bounds of $T$ is the least upper bound of $T$.
\par
(c) Let $x'$ be the greatest lower bound of the co-faces of $x$.
As the set of co-faces of  $ x $ is bounded below by $x$, one has $x' \geq x$.
Since $S$ is manifold--like, a non-maximal cell $x$ has at-least two distinct
co-faces $z_1$ and  $z_2$. But then $x' \leq z_1 \wedge z_2$, implying
$\rk(x') \leq \rk(z_1 \wedge z_2) < \rk(z_i) = \rk(x) + 1$.
Hence $\rk(x') \leq \rk(z_1 \wedge z_2) \leq \rk(x)$. It follows that $x' = x = z_1 \wedge z_2$.
\par
(d) Use induction on $(\rk(y) - \rk(x))$. Axiom (3) implies that
any cell of rank at-least 1 has at-least two faces, which proves part (d), for $\rk(y) - \rk(x) = 1$.
Suppose $\rk(y) - \rk(x) = k$ and assume the result for all $x< y$ with $\rk(y) - \rk(x) < k$.
By the induction hypothesis, $y$ has a facet $z$ of co-dimension 2, such that $z \notin U(x)$.
Of the two cells in between $y$ and $z$, at-least one must not be above $x$, thus providing us
with a face of $y$, that does not belong to $U(x)$.
\end{proof}
%
%
\begin{definition/lemma}
\label{Sop}
Let $S$ be a combinatorial cell complex. Assume $S$ is manifold--like, of dimension $n$.
For each $x \in S$, introduce a new symbol $x^{\circ}$, to be called the {\it dual cell} of $x$.
Let $S^{\circ} = \lbrace x^{\circ} : x \in S \rbrace$ with the partial order defined by
$x^{\circ} \leq^{\circ} y^{\circ} $ if and only if $ x \geq y$. Define a rank
function on $S^{\circ}$ by $\rk^{\circ}(x^{\circ}) = n - \rk(x)$.
It follows from lemma \ref{firstlemma} that $S^{\circ}$ is a combinatorial cell complex. 
It is called the {\it dual c.c.c} of $S$. The $r$--cells of $S^{\circ}$ correspond to 
the $(n-r)$--cells of $S$. 
The non-singularity of $S$ implies that $S^{\circ}$ is non-singular.
The boundary and co-boundary of $S$ are respectively the co-boundary and boundary of $S^{\circ}$.
Thus, if $S$ is manifold--like, then $S^{\circ}$ is also manifold--like and $(S^{\circ})^{\circ} = S$.
\end{definition/lemma}
%
%
\begin{remark}
From lemma \ref{firstlemma}(a), we see
in particular, that every cell is the least upper bound of its vertices. So we can identify
each cell with its set of vertices.
Thus, to define a c.c.c, we can start from the vertex set $S_0$, specify the subsets
of $S_0$ which correspond to the cells and the rank of each cell. The partial order is induced
by inclusion. It will be sometimes convenient to think of the empty set
$\emptyset$ as a cell of rank $-1$, lying below every vertex and consider the partially ordered
set $\tilde{S} = S \cup \lbrace \emptyset \rbrace $. Of-course $ \tilde{S} $ is not a
c.c.c.
\end{remark}
%
%
%
%
\section{Orientation on a combinatorial cell complex}
\begin{figure}[ht!]
\begin{center}
\includegraphics{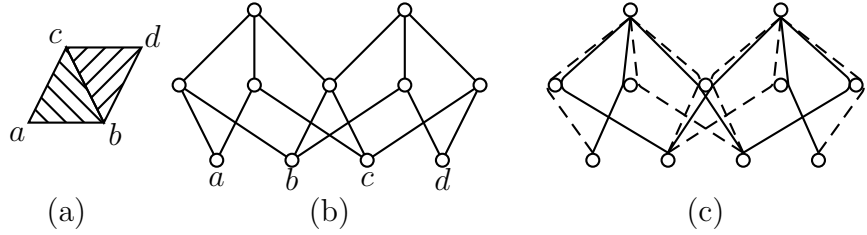}
\caption{Example of a 2-dimensional c.c.c: (a) shows two triangles joined along a common edge.
(b) shows the partially ordered set of the c.c.c corresponding to this geometric figure. 
(c) shows the flags of the c.c.c, drawn in two kinds of lines, showing an orientation.}
\label{twotriangle}
\end{center}
\end{figure}

\begin{figure}[ht!]
\begin{center}
\includegraphics{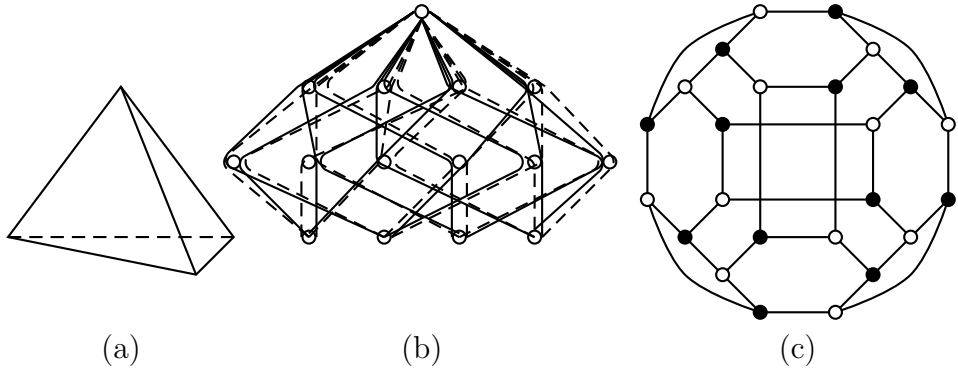}
\caption{a three dimensional c.c.c: (a) the tetrahedron. 
(b) the flags drawn in two kind of lines. (c) the adjacency graph of flags.}
\label{tetrahedron}
\end{center}
\end{figure}

\begin{figure}[ht!]
\begin{center}
\includegraphics{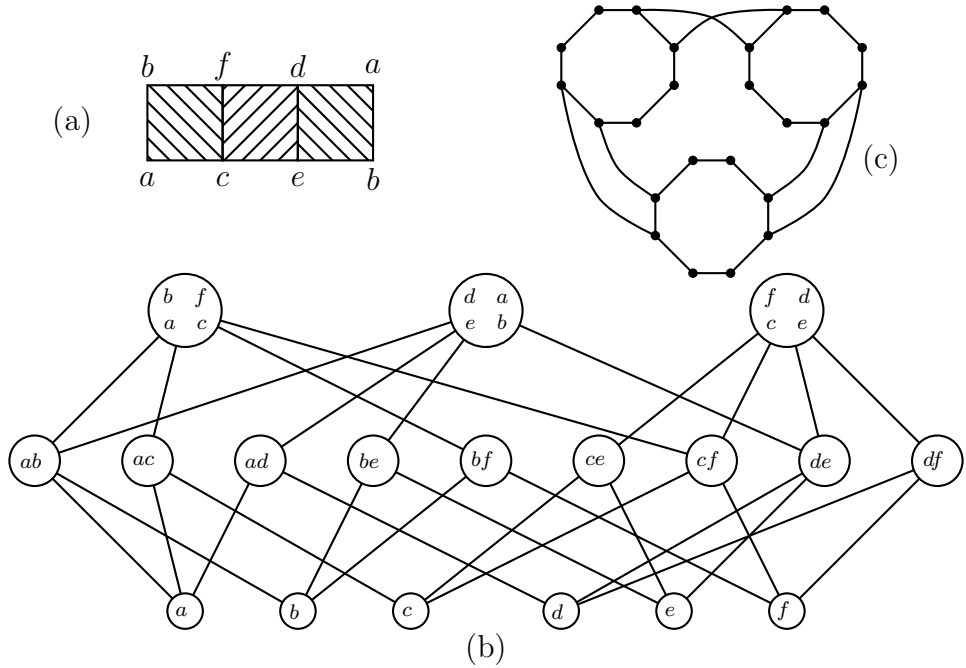}
\caption{(a) shows the Mobius strip broken up into three squares.
(b) shows the c.c.c corresponding to the Mobius strip. (c) shows the adjacency graph of flags;
this graph is not bipartite.}
\label{mobius}
\end{center}
\end{figure}

\begin{figure}[ht!]
\begin{center}
\includegraphics{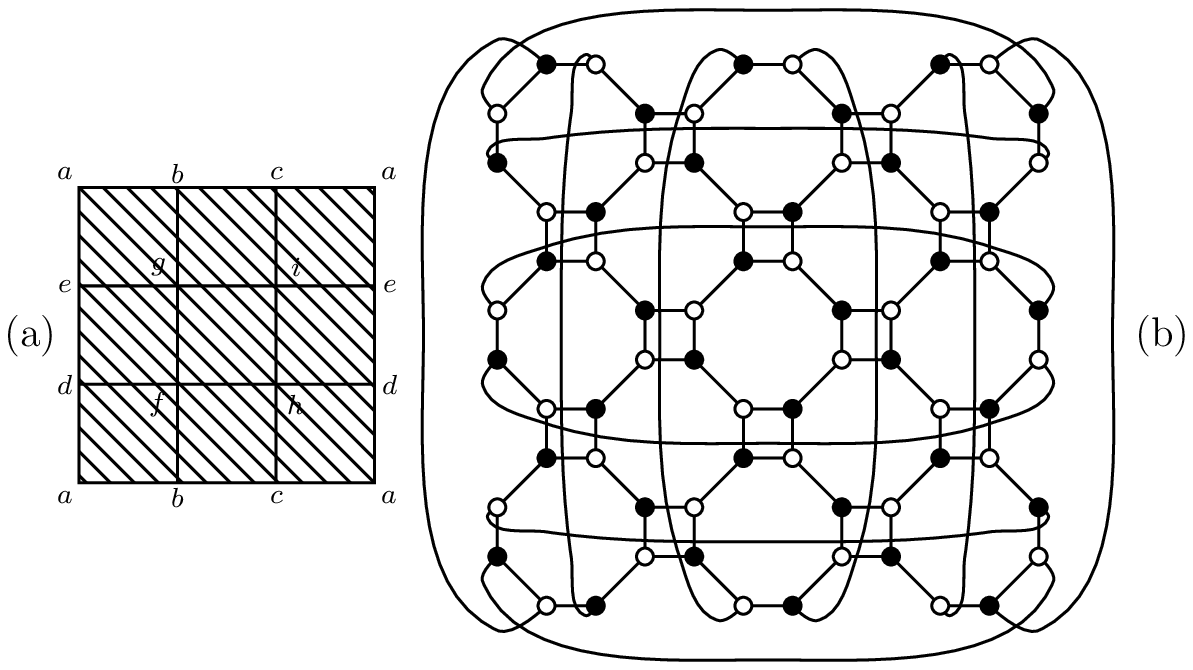}
\caption{(a) the torus broken up into 9 squares. (b) the adjacency graph of flags.}
\label{torus}
\end{center}
\end{figure}

%
%
%
%
\begin{definition}
\label{orientationdef}
Let $S$ be an equidimensional c.c.c, of dimension $n$. 
In particular, $S$ is a poset. So one has the usual notion of 
the barycentric subdivision of $S$. The (first) {\it barycentric subdivision} of $S$,
denoted by $S^{(1)}$, is the set of all totally ordered subsets of $S$. 
The barycentric subdivision of $S$, with partial order induced by inclusion,
is a c.c.c (in-fact a simplicial complex). The $r$--cells of $S^{(1)}$ are 
\begin{equation*}
S^{(1)}(r) = \lbrace \lbrace x_0 < x_1 < \dotsb < x_r\rbrace \colon x_j \in X \rbrace.
\end{equation*}
A {\it flag} in $S$ is an $n$--cell of $S^{(1)}$. In other words, a flag in $S$ is
a maximal totally ordered subset $\lbrace x_0 < x_1 < \dotsb < x_n \rbrace$ of $S$
such that $x_i \in S(i)$. Let $\cF(S)$ be the set of flags in $S$. 
We use the abbreviations $\cF(x) = \cF(\cl_S(x))$ and $\cF(x^{\circ}) = \cF(\cl_{S^{\circ}}(x^{\circ}))$.
A flag in $\cF(x)$ is called a {\it flag below} $x$. A flag in $\cF(x^{\circ})$ is called a {\it flag above} $x$.
\par
Two flags $F_1$ and $F_2$ are called {\it adjacent} if they differ only in one step,
that is, if the corresponding $n$--cells of $S^{(1)}$ have a common face.
The adjacency graph\footnote{a graph is a one dimensional CW-complex.} of flags
in $S$ will also be denoted by $\cF(S)$.
The vertices of this graph are the flags in $S$. Two flags are joined by an
edge if and only if the two flags are adjacent. 
\par
We say that $S$ is {\it flag--connected} if $\cF(S)$ is a connected graph.
We say that $S$ is {\it orientable} if the graph $\cF(S)$ is connected and bipartite.
An {\it orientation} $\omega$ on $S$ is a coloring of the flags in $S$ with two colors such that
adjacent flags get opposite color.
In other words, an orientation $\omega$ on $S$ is a function $\omega \co \cF(S) \to \lbrace \pm 1 \rbrace$,
such that $\omega(\gamma) = - \omega(\gamma')$ if $\gamma$ and $\gamma'$ are adjacent flags.
Since the graph $\cF(S)$ is assumed to be connected, an orientable c.c.c $S$ has
two possible orientations.
\par
Let $x \in S$. If $\cl (x)$ is flag--connected (resp. orientable), we say that $x$ is flag--connected
(resp. orientable). An orientation on $\cl (x)$ is referred to as an orientation on $x$.
\end{definition}
%
%
\begin{example} The above definition of orientation is central to our work. So we pause to 
illustrate the definition through examples of a few non-singular c.c.c's,
shown in the figures \ref{twotriangle}, \ref{tetrahedron}, \ref{mobius} and \ref{torus}.
The flags that map to $1$ are drawn in solid lines or solid dots and the ones that map to $-1$ are
drawn in dotted lines or hollow dots. Interchanging the solid lines (resp. solid dots) and the
dotted lines (resp. hollow dots), one gets the reverse orientation.
\end{example}
%
%
%
%
%
%
%
%
%
%
%
\begin{remark}
\label{orientablecellremark}
\begin{enumerate}
\item
Suppose $S$ is a c.c.c with flag--connected cells. Suppose $x$ is a cell of $S$ and
$y$ is a face of $x$. Then each flag below $y$ can be extended uniquely to a flag below $x$.
So an orientation $\omega$ on $x$ induces an orientation $\omega\vert_y$ on $y$, defined by
\begin{equation*}
 \omega\vert_y (\gamma) = \omega(\gamma \cup \lbrace x \rbrace) \text{\;\; for \;\;} \gamma \in \cF(y).
\end{equation*}
It follows that, if each maximal cell of $S$ is orientable, then each cell of $S$ is orientable.
\par
An orientation on $S$ determines an orientation on each maximal cell of $S$. So if $S$ is
orientable, with flag--connected cells, then each cell of $S$ is orientable.
\item If two cells $x_+$ and $x_-$ share a face $x$, then an orientation on
$\cl \lbrace x_+, x_- \rbrace$ induces two opposite orientations on $ \cl (x)$,
one coming from the orientation on $\cl (x_+) $ and the other one coming from the orientation on
$\cl (x_-) $.
\item Notice that an orientable c.c.c must be non-singular. If an $1$--cell of $S$ has $r$ faces,
or if there are $r$ maximal cells of $S$ sharing a common face, then the
graph $\cF(S)$ contains a complete graph on $r$ vertices. So $\cF(S)$ can be
bipartite only if $r \leq 2$.
\item Suppose $S$ is a non-singular c.c.c with only one maximal cell.
Then an orientation on $S$ determines an orientation on the boundary of $S$.
\end{enumerate}
\end{remark}
%
%
\begin{definition}
\label{defsigns}
 Let $x$ be an orientable cell of a c.c.c $S$ and $y$ be an orientable face of $x$.
Let $\omega$ be an orientation on $x$ and $\mu$ be an orientation on $y$.
We define,
\begin{equation*}
 s(\omega, \mu) = \begin{cases} 1 & \text{\;\;if\;\;} \omega\vert_y = \mu, \\
                   -1 & \text{\;\;if\;\;} \omega\vert_y = - \mu.
                  \end{cases}
\end{equation*}
If $\omega_x$ is an orientation on $x$ and $\omega_y$ is an orientation on $y$, then we
write $s(x,y) = s(\omega_x, \omega_y)$. To determine $s(x,y)$,
consider a flag $\gamma \in \cF(x)$ of the form $\gamma = \lbrace x > y > \dotsb \rbrace$. Then
\begin{equation}
s(x,y) = \omega_x(\gamma)/\omega_y(\gamma \setminus \lbrace x \rbrace).
\label{sbyflag}
\end{equation}
Since the graph $\cF(y)$ is connected, the right hand side
of equation \eqref{sbyflag} does not depend on the choice of the flag $\gamma$.
\end{definition}
%
%
%
%
\section{homology and cohomology groups}
%
%
%
%
\label{sechomcohom}
\begin{topic}
For this section, let $S$ be a c.c.c such that each cell of $S$ is orientable.
Pick an orientation on each cell $x$ of $S$, denoted by
$\omega_x \co \cF( \cl (x)) \to \lbrace \pm 1 \rbrace$.
Given this data, we can associate a sign $s(x,y) \in \lbrace \pm 1 \rbrace$,
for each pair $x$ and  $y$, where $x$ is a cell and $y$ is a face of $x$ (see \ref{defsigns}).
The key equation satisfied by the numbers $s(x,y)$ is given in the following lemma. 
Axiom (4) in the definition of a c.c.c, which is our main axiom, is used here.
\end{topic}
%
%
\begin{lemma}Given the setup in section \ref{sechomcohom} so far,
Let $z$ be a co-dimension $2$ facet of $x \in S$. Let
 $y_+$ and $y_-$ be the two cells in between $x$ and $z$, that is,
$\Delta x \cap \nabla z = \lbrace y_+ , y_- \rbrace$.
Then
\begin{equation}
s(x,y_+)  s(y_+, z) + s(x,y_-) s(y_-, z) = 0.
\end{equation}
\label{rhombus}
\end{lemma}
\begin{proof}
Let $\gamma = \lbrace z = z_0 > z_1 > \dotsb \rbrace$ be a flag below $z$. Let
$\gamma_+ = \lbrace x > y_+ > z_0 > z_1 > \dotsb \rbrace$
and $\gamma_- = \lbrace x > y_- > z_0 > z_1 > \dotsb \rbrace$ be the two flags below $x$
that extend $\gamma$. Then
\begin{equation*}
s(x,y_+) s( y_+, z) =
\frac{\omega_x (\gamma_+)}{\omega_{y_+}(\gamma_+ \setminus \lbrace x \rbrace ) } \cdot
\frac{\omega_{y_+}(\gamma_+ \setminus \lbrace x \rbrace ) }{\omega_z(\gamma)}
= \frac{\omega_x(\gamma_+)}{\omega_z(\gamma)}
\end{equation*}
Similarly $s(x,y_-) s( y_-, z) =  \omega_x(\gamma_-)/\omega_z(\gamma)$.
Since $\gamma_+$ and  $\gamma_-$ are adjacent flags in $\cF(x)$,
the lemma follows.
\end{proof}
%
%
\begin{definition}
Now we can define chain complexes, boundary maps, homology groups et-cetera in the standard
fashion.
For each cell $x$ of $S$, we introduce a formal variable, denoted by $[x]$.
The group of $i$--chains in $S$ with integer coefficients, denoted by $C_i(S)$, is the free
$\ZZ$--module with basis $\lbrace [x] \colon x \in S(i) \rbrace$.
(Of course, one can replace $\ZZ$ by other commutative rings but we shall restrict ourselves to
integer coefficients). Let
\begin{equation*}
\partial [x] = \sum_{ y \in \Delta x } s(x,y) [y]
\text{\;\;\;\; and \;\;\;\;}
\delta [x]  = \sum_{ z \in \nabla x } s(z,x) [z].
\end{equation*}
Define the {\it boundary map} $\partial \co C_i(S) \to C_{i-1}(S)$ and the {\it co-boundary} map
$\delta \co C_i(S) \to C_{i+1}(S)$ by linearly extending the above.
In other words, for an $i$--chain $\sigma = \sum_{x \in S(i) } r_x [x]$, let
\begin{equation*}
\partial (\sum_{x \in S(i) } r_x [x]) = \sum_{x \in S(i) } r_x \partial [x]
\text{\;\; and \;\;}
\delta (\sum_{x \in S(i) } r_x [x]) = \sum_{x \in S(i) } r_x \delta [x].
\end{equation*}
The image of a minimal (resp. maximal) cell under the boundary (resp. co-boundary) map
is defined to be zero.
If $\sigma \in C_i(S)$ such that $\partial \sigma = 0$ (resp. $\delta \sigma = 0$) we say that
$\sigma$ is an $i$--cycle (resp. $i$--cocycle).
\end{definition}
%
%
\begin{lemma} Given the setup in section \ref{sechomcohom} so far, one has $\partial^2 = 0$ and $\delta^2 = 0$. \end{lemma}
\begin{proof}
The proof follows from axiom (4) in the definition \ref{cccdef} and lemma \ref{rhombus}.
\end{proof}
%
%
\begin{definition} 
\label{defhom}
Let $C_i = C_i(S)$. The lemma above shows that
$(C_i, \partial)$ and $(C_i, \delta)$ are chain complexes.
We define the cellular homology  (resp. cellular cohomology) of
$S$ to be the homology of the chain complex $(C_i, \partial)$,
(resp. $(C_i, \delta)$).
\begin{equation*}
H_i(S)  = \frac{\op{ker}(\partial \co C_i \to C_{i-1})}{\op{im}(\partial \co C_{i+1} \to C_i)}
\text{\;\; and \;\;}
H^i(S)  = \frac{\op{ker}(\delta \co C_i \to C_{i+1})}{\op{im}(\delta \co C_{i-1} \to C_i)}.
\end{equation*}
\end{definition}
%
%
\begin{remark}$\;\;$
\begin{enumerate}
\item To define the homology and cohomology of $S$, we need each cell of $S$
to be orientable. We do not require that $S$ is non-singular or even equidimensional.
If each cell of $S$ is orientable, and $T$ is a closed subset of $S$, then each cell of $T$ is
also orientable. So the homology/cohomology groups of $T$ are well defined. However $T$ need not
be equidimensional or non-singular, even if $S$ were. We  shall have
occasion to consider homology groups of such $T$.
\item Suppose $S$ is a c.c.c with orientable cells.
Given an orientation $\omega_y$ on each cell $y$ of $S$,
we get the chain complex $(C_{\bullet}, \partial)$ as defined above. Let us temporarily
write $(C_{\bullet}, \partial) = (C_{\bullet}^{\omega}, \partial^\omega)$ to emphasize that
the chain complex depends on the choice of $\omega_y$'s. However,
as we shall now see, choosing a different set of orientations,
give an isomorphic chain complex. Let $\lbrace \mu_y \colon y \in S \rbrace $ be another
set of orientations on the cells of $S$. Define $t(y) = 1$ if $\omega_y = \mu_y$
and $t(y) = -1$ if $\omega_y = - \mu_y$.
Then it can be easily checked that the map $[y] \mapsto t(y) [y]$ gives an isomorphism, 
\begin{equation*}
 (C_{\bullet}^{\omega}, \partial^\omega) \simeq  (C_{\bullet}^{\mu},
 \partial^\mu),
\end{equation*}
of chain complexes.
So the homology groups do not depend on the choice of $\omega_y$. The same remark applies
to the cohomology groups.
\item Assume that $S$ has orientable cells. Then each $1$--cell has two vertices.
The zero cycles of $S$ are just linear combinations of vertices of $S$.
Usually we shall assume that $\omega_v(\lbrace v \rbrace) = 1$
for each cell $v$ of rank zero. Under this assumption, if $v_+$ and $v_-$ are the two vertices of
an $1$--cell $x$, then $s(x, v_+) + s(x, v_-) = 0$. So two vertices $v_1$ and $v_2$ are in
the same homology class if and only if they can be ``joined by a sequence of $1$--cells''.
\par
Consider the graph $S_{\leq 1}$ whose edges correspond to the $1$--cells of $S$
and the two endpoints of an edge $x$ correspond to the two rank zero cells of $S$ below $x$.
Then $H_0(S)$ is simply the zero-th homology of the  one dimensional
CW--complex $S_{\leq 1}$. 
Suppose the graph $S_{\leq 1}$ has $r$ connected components. 
Then $H_0(S)$ is a free abelian group of rank $r$.
If one vertex is chosen from each component of the graph $S_{\leq 1}$, then $H_0(S)$ is
freely generated by the homology classes of these $r$ vertices. In particular,
if $H_0(S) \simeq \ZZ$, then $H_0(S)$ is generated by the class of any vertex of $S$.
\item Let $T$ be a closed subset of $S$. Let $C_i(S,T) = C_i(S)/C_i(T)$. If $\sigma \in C_i(T)$,
then its boundary $\partial \sigma$ belongs to $C_{i-1}(T)$. Thus $\partial$ induces boundary
maps $\partial_T^S \co C_i(S, T) \to C_{i-1}(S,T) $.
We define the {\it relative homology} of the pair $(S,T)$ to be homology of the chain complex
$(C_i(S, T), \partial_T^S)$.
\end{enumerate}
\end{remark}
%
%
\begin{lemma}
\label{simphom}
Let $S$ be a simplicial complex. For each simplex $\gamma = \lbrace x_0, \dotsb, x_r \rbrace \in S$ of rank $r$,
choose a total ordering, $x_r  <_{\gamma} x_{r-1} <_{\gamma} \dotsb <_{\gamma} x_0$,
on the set of vertices of $\gamma$. Assume that these total orderings are compatible with each other,
that is, if $\gamma' \subseteq \gamma$, then $<_{\gamma'}$ is the restriction of $<_{\gamma}$ to the
vertices of $\gamma'$.\footnote{For example, a total ordering on all the vertices of $S$,
induces a compatible family of total orderings on the vertices of each simplex of $S$.}
Now consider $S$ as a combinatorial cell complex. 
\par 
Then each cell of $S$ is flag--connected
and there exists an orientation $\omega_{\gamma}$ on each cell $\gamma$ of $S$
such that 
\begin{equation*}
s(\gamma, \gamma \setminus \lbrace x_i \rbrace) = (-1)^i.
\end{equation*}
It follows that the homology of the c.c.c $S$ (as defined in \ref{defhom})
coincides with the simplicial homology of the simplicial complex $S$
(as defined, for example, in section 3.2 of \cite{DK:CA}).
\end{lemma}
\begin{proof}
Let $\gamma = \lbrace x_0, \dotsb, x_r \rbrace$ be a simplex of $S$. A total ordering,
given by $x_r <_{\gamma} \dotsb <_{\gamma} x_0$, on the vertices of $\gamma$,
induces an orientation on $\gamma$, as follows.
\par
Given a flag $\Gamma = \lbrace \gamma = \Gamma_0 \supset \Gamma_1 \supset \dotsb \supset \Gamma_r \rbrace$
in $\cl(\gamma)$, one gets a permutation $P_{\gamma}(\Gamma)$ of $(r+1)$ letters,
defined by $\Gamma_i \setminus \Gamma_{i+1} = \lbrace x_{P_{\gamma}(\Gamma)(i)} \rbrace$.
Define $\omega_{\gamma} \co \cF(\gamma) \to \lbrace \pm 1 \rbrace$ by
\begin{equation*}
\omega_{\gamma}(\Gamma) = \op{sign}(P_{\gamma}(\Gamma)).
\end{equation*}
If $\Gamma$ and $\Gamma'$ are adjacent flags below $\gamma$, then
the permutations $P_{\gamma}(\Gamma)$ and $P_{\gamma}(\Gamma')$ differ
by a transposition. So $\omega_{\gamma}$ is an orientation on the cell $\gamma$.
Notice that $\cF(\gamma)$ is flag connected since the symmetric group is
generated by transpositions.
\par
 To determine $s( \gamma, \gamma \setminus \lbrace x_0 \rbrace)$, 
consider the flag $\Gamma = \lbrace \Gamma_0 \supseteq \Gamma_1 \supseteq \dotsb \supseteq \Gamma_r \rbrace$
given by $\Gamma_i = \gamma \setminus \lbrace x_0, \dotsb, x_{i-1} \rbrace$.
Then $P_{\gamma}(\Gamma)$ and
$P_{\gamma \setminus \lbrace x_0 \rbrace } ( \Gamma \setminus \lbrace \gamma \rbrace )$
are both equal to the identity permutation. So
\begin{equation*}
\omega_{\gamma}(\Gamma)
= \omega_{\gamma \setminus \lbrace x_0 \rbrace } ( \Gamma \setminus \lbrace \gamma \rbrace ) = 1.
\end{equation*}
It follows that $s( \gamma, \gamma \setminus \lbrace x_0 \rbrace)
= \omega_{\gamma}(\Gamma)/\omega_{\gamma \setminus \lbrace x_0 \rbrace } ( \Gamma \setminus \lbrace \gamma \rbrace ) = 1$.
To compare $s( \gamma, \gamma \setminus \lbrace x_i \rbrace)$ and $s( \gamma, \gamma \setminus \lbrace x_{i+1} \rbrace)$, consider two adjacent flags $\Gamma_+$ and $\Gamma_-$ in $\cF (\gamma )$, having the following form:
\begin{align*}
\Gamma_+ &= \lbrace \gamma \supseteq \gamma \setminus \lbrace x_i \rbrace \supseteq \gamma \setminus \lbrace x_i , x_{i+1} \rbrace \supseteq \dotsb \rbrace, \\
\Gamma_- &= \lbrace \gamma \supseteq \gamma \setminus \lbrace x_{i+1} \rbrace \supseteq \gamma \setminus \lbrace x_i , x_{i+1} \rbrace \supseteq \dotsb \rbrace.
\end{align*}
Since $\Gamma_+$ and $\Gamma_-$ are adjacent flags, we have $\omega_{\gamma}(\Gamma_+) = - \omega_{\gamma}(\Gamma_-)$.
On the other hand, the flags 
$\Gamma_+ \setminus \lbrace \gamma \rbrace \in \cF( \gamma \setminus \lbrace x_i \rbrace)$
and $\Gamma_- \setminus \lbrace \gamma \rbrace \in \cF( \gamma \setminus \lbrace x_{i+1} \rbrace)$
correspond to the same permutation. Hence
$\omega_{\gamma \setminus \lbrace x_i \rbrace} ( \Gamma_+ \setminus \lbrace \gamma \rbrace ) =
\omega_{\gamma \setminus \lbrace x_{i+1} \rbrace} ( \Gamma_- \setminus \lbrace \gamma \rbrace )$.
It follows that $s( \gamma, \gamma \setminus \lbrace x_i \rbrace)$
and $s( \gamma, \gamma \setminus \lbrace x_{i+1} \rbrace)$ have opposite signs.
\end{proof}
%
%
%
%
\section{Stellar subdivision}
\begin{figure}[ht!]
\begin{center}
\includegraphics{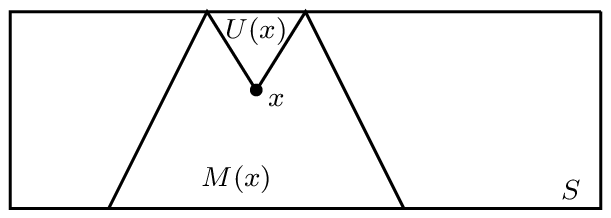}
\caption{}
\label{M}
\end{center}
\end{figure}

\begin{figure}[ht!]
\begin{center}
\includegraphics{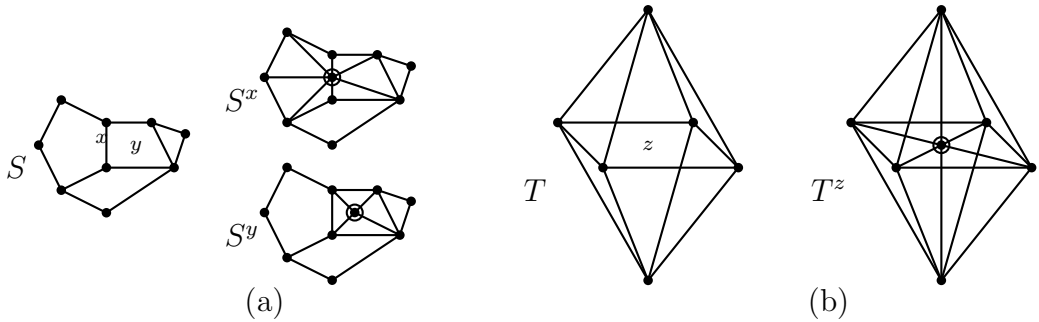}
\caption{ (a) shows a two dimensional c.c.c $S$, the stellar subdivision 
$S^x$ at the $1$--cell $x$ joining the square and the pentagon
and the stellar subdivision $S^y$ at the 2-cell $y$.
(b) shows a three dimensional c.c.c $T$ and it subdivision $T^z$, where
the 2-cell $z$ is the square in the middle.
The new vertices, namely $C_x(\emptyset)$, $C_y(\emptyset)$ and $C_z(\emptyset)$,
are marked with a circle.}
\label{subdiv}
\end{center}
\end{figure}

%
%
%
%
\label{secsteller}
We would like to show that if $S$ is a manifold--like c.c.c with orientable and acyclic
cells, then the homology of $S$ is isomorphic to that of its barycentric subdivision $S^{(1)}$.
It is easy to write down a chain map from the $i$--chains of $S$ to those of $S^{(1)}$.
But it seems difficult to show directly that this map induces isomorphism of homology groups, 
since the cell structure of $S^{(1)}$ is very different from the cell structure of $S$.
For this purpose, we want to break up the transition from $S$ to $S^{(1)}$ into many successive
``stellar subdivisions" or ``stellar refinements". 
In each step, the cell structure is only ``locally" modified. This makes it easier to compare the
homology groups in successive steps. 
Stellar subdivisions of simplicial and cell complexes arise in many places in literature,
for example, see \cite{ES:SS}, \cite{DK:CA}.
\begin{definition}
Let $x$ be a cell of a c.c.c $S$. Define the {\it star} of $x$ to be
\begin{equation*}
\st(x) = \cl( U(x)).
\end{equation*}
We also define $M(x) = \st(x) \setminus U(x)$ and $\tilde{M}(x) = M(x) \cup \lbrace \emptyset \rbrace$
(see figure \ref{M}).
Both $\st(x)$ and $M(x)$ are closed subsets of $S$. So these are sub--c.c.c's of $S$. When there is a
possibility of confusion, we write $\st_S(x)$ and $M_S(x)$.
Say that $S$ is a {\it star around $x$}, if $\st_S(x) = S$.
\end{definition}
\begin{definition}
\label{stellerdef}
Let $S$ be a c.c.c and $x \in S(i)$ for some $ i \geq 1$.
We want to define a new c.c.c $S^x$, to be called {\it the stellar subdivision of $S$ at $x$}.
(To get the idea, look at the examples in figure \ref{subdiv}).
For each $ y \in \tilde{M}(x) $, introduce new cells $C_x(y)$, to be called
{\it the cone over $y$ with vertex at $x$}. Define
$S^x(0) = S(0) \cup \lbrace C_x(\emptyset) \rbrace$ and
\begin{equation*}
S^x(r) = \lbrace y \in S(r) \colon y \ngeq x  \rbrace \cup
\lbrace C_x(y) : y \in \tilde{M}(x)(r -1)  \rbrace ,
\end{equation*}
with the convention that $\tilde{M}(x)(-1) = \lbrace \emptyset \rbrace$.
There are two kinds of cells in $S^x$.
The first kind consists of the cells of $S \setminus U_S(x)$; these will be called
the {\it old cells}.
The second kind consists of the cones; these will be called the {\it new cells}.
\par
Next, we define the partial order on $S^x$.
Given two cells $y$ and $z$ of $S^x$, the relation $y \leq_{S^x} z$ holds if and only if one of
the following conditions hold. 
\begin{itemize}
\item Both $y$ and $z$ are old cells and $y \leq_S z$.
\item $y$ is an old cell, $z = C_x(z')$ is a new cell and $ y \leq_S z'$.
\item Both $y = C_x(y')$ and $z = C_x(z')$ are new cells and $y' \leq_S z'$.
\end{itemize}
We shall check in a moment that $S^x$ is a c.c.c.
If $T$ is obtained from $S$ by successive stellar refinements, then we
say that $T$ is a refinement of $S$.
\end{definition}
\begin{remark}
\label{stellerremark}
\begin{enumerate}
\item 
Let $T = \st_S(x)$ and $M(x) = M_S(x) = M_T(x)$.
There is a canonical isomorphism: $\st_{S}(x)^x \simeq \st_{S^x}(C_x(\emptyset))$.
On both sides, the $r$--cells are 
\begin{equation*}
M(x)(r) \cup \lbrace C_x(y) \colon y \in \tilde{M}(x)(r-1) \rbrace.
\end{equation*}
On both sides, the partial order and rank are defined in the same way.
We shall often identify $\st_S(x)^x$ as a sub-c.c.c of $S^x$, via the above isomorphism.
\item Taking a stellar refinement at $x$ only changes the cell structure ``around $x$''.
More precisely, $\st_S(x) \subseteq S$ is replaced by $\st_S(x)^x \simeq \st(C_x(\emptyset)) \subseteq S^x$.
The rest of the cell structure remains unchanged.
\item The cells in $U(x) \subseteq S$ ``die" in the process of
stellar subdivision at $x$.
The rest of the cells of $S$ ``survive" as cells of $S^x$; these are the
old cells. Finally, for each cell  $y \in \tilde{M}(x)$, a cell called $C_x(y)$
is ``born"; these are the new cells. For later use, we note the following.
\begin{itemize}
\item There are no new cells below an old cell.
\item Among the faces of $C_x(y)$, there is only one old cell, namely $y$ itself.
\end{itemize}
\item 
While defining $S^x$, we have assumed that the rank of $x$ is atleast one, because
this is the only case we shall need.
However, the definition makes sense even when $x$ is a cell of rank zero.
In this case the vertex $x$ gets replaced by the vertex $C_x(\emptyset)$.
\end{enumerate}
\end{remark}
\begin{lemma}
Let $S$ be a c.c.c and $x$ be a cell of $S$ of rank at-least one. Then,
\newline
(a) $S^x$ is a c.c.c.
\newline
(b) If $S$ is equidimensional, of dimension $n$, then so is $S^x$.
\newline
(c) If each $1$--cell of $S$ has two vertices, then the same is true for
each $1$--cell of $S^x$.
\newline
(d) Suppose $S$ is equidimensional, of dimension $n$. If there are at-most two (resp. exactly two)
$n$--cells above each $(n-1)$--cell of $S$, then the same is true for $S^x$.
\newline
(e) If $S$ is non-singular (resp. manifold--like), then $S^x$ is non-singular (resp. manifold--like).
\label{stellerfirstlemma}
\end{lemma}
The proof, given in appendix \ref{stellerfirstlemmapf}, is easy but a little tedious.
It is mainly because we have to separate the argument into cases, depending on
whether the cell of $S^x$ we are dealing with is a cone or not.
\par
%
%
%
We shall have occasion to consider repeated stellar subdivision of a c.c.c.
We shall write $(X^x)^y = X^{x y}$.
The c.c.c one obtains by repeated stellar subdivision depends,
in general, on the order in which the subdivision points are chosen.
However, we have the following result.
\begin{lemma}
\label{commute}
Let $X$ be a c.c.c and $\lbrace x_1 , \dotsb, x_k \rbrace \subseteq X$ such that
$U_X(x_i) \cap U_X(x_j) = \emptyset$ for all $ i \neq j$. Then the refinement
$X_{(k)}  = X^{x_1  x_2 \dotsb x_k }$ has the following description:
\begin{equation*}
X_{(k)} = \op{\cup}_{j=1}^k  \lbrace C_{x_j}(v) \colon v \in \tilde{M}_X(x_j) \rbrace \cup
\bigl( X \setminus \op{\cup}_{ j = 1}^k U_X(x_j) \bigr).
\end{equation*}
As before, the cells of the form $C_{x_j}(v)$ are called the new cells and the rest are called the old cells.
The partial order on $X_{(k)}$ is defined by the following rules. One has $\alpha \leq_{X_{(k)}} \beta$ if
and only if one of the following three conditions hold:
\begin{itemize}
\item both $\alpha$ and $\beta$ are old and $\alpha \leq_X \beta$.
\item $\alpha$ is old, $\beta = C_{x_j}(\beta')$ is new and $\alpha \leq_X \beta'$.
\item both $\alpha$ and $\beta$ are new, there is a $j$ between $1$ and $k$ such that
$\alpha = C_{x_j}(\alpha')$, $\beta = C_{x_j}(\beta')$ and $\alpha' \leq_X \beta'$.
\end{itemize}
It follows from this description that there are no old cells above a new cell and
$X_{(k)}$ does not depend on the order of subdivision.
\end{lemma}
The proof is given in appendix \ref{commutepf}.
\par
%
%
Suppose $X$ is a c.c.c such that each cell of $X$ is orientable but $X$ itself is not orientable.
We will need to consider the homology groups of such an $X$ and of its stellar subdivision $X^x$.
We need the following lemma to make sure that the homology of $X^x$ is well defined.
\begin{lemma}
\label{cellstayorientable}
Let $X$ be a c.c.c and let $x$ be a cell of $X$. 
\par
(a) If each cell of $X$ is flag-connected, then each cell of $X^x$ is flag connected.
\par
(b) If each the cell of $X$ is orientable, then each cell of $X^x$ is also orientable.
More precisely, one has the following:
Let $y \in M_X(x)$ with $\rk(y) = n-1$.
Let $S = \cl_X(y) \subseteq X$ and $S' = \cl(C_x(y)) \subseteq X^x$.
Given a flag $\gamma \in \cF(S')$, there is an $i \geq 0$ such that
\begin{equation*}
\gamma = \lbrace C_x(y_0) > C_x(y_1) > \dotsb  > C_x(y_i) > y_i  > y_{i+2} > y_{i+3} >  \dotsb > y_n \rbrace
\end{equation*}
where $y_0 = y$ and $y_j \in S$, with the exception that $y_n = \emptyset$ if $i = n$.
We let $l(\gamma) = i$ and
\begin{equation*}
\tilde{\gamma} = \lbrace y_0 > y_1 > \dotsb > y_i  > y_{i+2} > y_{i+3} > \dotsb > y_n \rbrace \in \cF(S),
\end{equation*}
with the convention that $y_n = \emptyset$ is omitted if $i = n$.
If $\omega_y$ is an orientation on $S = \cl_X(y)$, then $\omega_{S'}$, defined by
\begin{equation*}
\omega_{S'}(\gamma) = (-1)^{l(\gamma)} \omega_y(\tilde{\gamma}),
\end{equation*}
is an orientation on $S' = \cl(C_x(y))$.
\end{lemma}
The proof is given in appendix \ref{cellstayorientablepf}.
%
%

\begin{definition}
Let $S$ be a c.c.c with orientable cells and $x \in S$. 
Fix an orientation $\omega_z$ for each cell $z \in S$. Given this data,
we define an orientation on each cell of $S^x$ as follows.
If $z \in S^x$ is an old cell, then $\cF_{S}(z) = \cF_{S^x}(z)$. So $\omega_z$ is already defined.
If $C_x(y)$ is a cone in $S^x$, then choose $\omega_{C_x(y)}$ as prescribed by lemma
\ref{cellstayorientable}(b). For a flag $\gamma$ with top two cells $C_x(y)$ and $y$, we have,
in the notation of lemma \ref{cellstayorientable}, $\tilde{\gamma} = \gamma \setminus C_x(y)$ and $l(\gamma) = 0$,
so $\omega_{C_x(y)}(\gamma ) = \omega_y( \gamma \setminus \lbrace C_x(y) \rbrace )$.
In other words, in the notation of \ref{defsigns}, we have
\begin{equation}
s(C_x(y), y) = 1.
\label{scone1}
\end{equation}
Suppose $y \in M_S(x)$ and $z$ is a face of $y$. So $z$ is a co-dimension 2 facet of  $C_x(y)$.
The two cells in between $C_x(y)$ and $z$ are  $C_x(z)$ and  $y$. From lemma \ref{rhombus},
one has,
\begin{equation*}
s(C_x(y), C_x(z)) s(C_x(z), z) = - s(C_x(y), y) s(y,z). 
\end{equation*}
Since $s(C_x(u), u) = 1$
for all $u$, it follows that
\begin{equation}
s(C_x(y), C_x(z)) = - s(y,z).
\label{scone}
\end{equation}
\end{definition}
\begin{lemma} Let $S$ be a c.c.c with orientable cells and $x \in S$.
For each $w \in S$, let $\Delta_1 w = \Delta w \setminus U(x)$ and
$\Delta_2 w = \Delta w \cap U(x)$.
Define
\begin{equation*}
\varphi([w] ) =
\begin{cases} \sum\limits_{y \in \Delta_1 w } s(w, y) [C_x(y)] & \text {if \;} w \in U(x), \\
[w] & \text{otherwise.}
\end{cases}
\end{equation*}
Then $\varphi$ defines a chain map $ (C_{\bullet}(S), \partial) \to (C_{\bullet}(S^x), \partial)$ and
hence induces an homomorphism $H_i(\varphi) \co H_i(S) \to H_i(S^x)$.
\label{StoSx}
\end{lemma}
\begin{proof}
Suppose $w \in U(x)$ and $y \in \Delta_1(w)$. 
Let $Z$ be the set of co-dimension $2$ facets of $w$, that are not greater than or equal to $x$.
From the description of partial order on $S^x$ and equations \eqref{scone1} and \eqref{scone}, we have, 
\begin{equation*}
\partial(C_x(y))=[y] - \sum_{z \in \Delta y} s(y,z) [C_x(z)].
\end{equation*}
It follows that
\begin{equation*}
\partial(\varphi[w])
= \sum_{y \in \Delta_1 w } s(w,y) \partial [C_x(y)] 
 = \sum_{y \in \Delta_1 w } s(w,y) [y]
- \sum_{ z \in Z } \Bigl[ \sum_{y \in \nabla z \cap \Delta_1 w } s(w,y) s(y, z) \Bigr] [C_x(z)].
\end{equation*}
In the second term of the final expression, we are summing over all pairs $(y,z)$  
such that $y\in \Delta w$, $z \in \Delta y$ and $y \notin U(x)$. 
So the set of $z$ that appear in the expression are in $Z$.
\par
Given $z \in Z$, let $y_+$ and $y_-$ be the two cells in between $w$ and $z$. Without
loss, we may assume that $y_+ \notin U(x)$.
We may write $Z$ as a disjoint union $Z = Z_1 \cup Z_2$, where $Z_1$ (resp. $Z_2$) consists of those $z \in Z$,
such that $y_- \in U(x)$ (resp. $y_- \notin U(x)$) (see figure \ref{Z_1}).
\begin{figure}
\begin{center}
\includegraphics{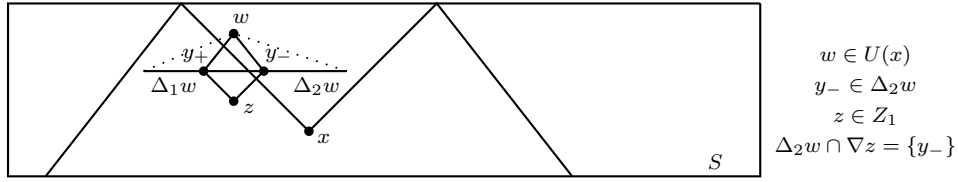}
\caption{the relevant cells around $z \in Z_1$}
\label{Z_1}
\end{center}
\end{figure}
For $z \in Z_2$, we have
$\sum_{y \in \nabla z \cap \Delta_1 w } s(w,y) s(y, z) = s(w,y_+)s(y_+, z) + s(w,y_-) s(y_-, z) = 0$.
It follows that
\begin{equation*}
\partial(\varphi[w]) =
\sum_{y \in \Delta_1 w } s(w,y) [y]
- \sum_{ z \in Z_1 } s(w,y_+) s(y_+, z) [C_x(z)].
\end{equation*}
To compute $\varphi(\partial[w])$, note that, if $y \neq y'$ are two cells in $\Delta_2 w$, then 
$\Delta_1 y \cap \Delta_1 y' = \emptyset$ and $\cup_{y \in \Delta_2 w} \Delta_1 y = Z_1$. It follows that
\begin{equation*}
\varphi(\partial[w])
 = \sum_{y \in \Delta_1 w \cup \Delta_2 w } s(w,y) \varphi [y] 
 = \sum_{y \in \Delta_1 w } s(w,y) [y] + \sum_{ z \in Z_1} s(w, y_-) s(y_-, z) [C_x(z)].
\end{equation*}
Using lemma \ref{rhombus} once more, we see that $ \partial \circ \varphi = \varphi \circ \partial$.
\end{proof}
%
%
%
%
\section{Lemmas on vanishing of homology groups}
%
%
%
%
\label{secvanishing}
\begin{definition}
\label{acyclicdef}
A c.c.c $S$ with orientable cells is {\it acyclic} if $H_i(S) = 0$ for $i > 0$
and $H_0(S) \simeq \ZZ$. As remarked in \ref{orientablecellremark}(3),
in such a situation, $H_0(S)$ is generated by the homology class of any vertex of $S$.
We say that $x$ is an acyclic cell if $\cl(x)$ is acyclic. In this section we want to
show that, if the cells of $S$ are acyclic, then the cells of $S^x$ are acyclic
and $H_{\bullet}(S) \simeq H_{\bullet}(S^x)$.
\end{definition}
\begin{lemma}
\label{cellstayacyclic}
Let $T$ be a c.c.c with orientable cells.
Let $x$ and $y$ be two cells of $T$ such that $y \in M(x)$.
Let $S = \cl_S(y)$ and $S' = \cl_{S^x}( C_x(y))$. Let
$j \co S \to S'$ be the inclusion map, $j(z) =z$. Then one has the following:
\newline
(a) The induced map on homology, $j_* \co H_i(S) \to H_i(S')$, is the zero map, for $i \geq 1$.
\newline
(b) If $y$ is acyclic, then so is $C_x(y)$.
\newline
(c) If all the cells of $T$ are acyclic, then all the cells of $T^{x}$ are also acyclic.
\end{lemma}
\begin{proof}
Let $z$ be a facet of $y$.
Since $x$ is not a facet of $y$, it is not a facet of $z$ either.
So $z$ remains a cell in $S^x$. So $j(z) = z$ defines an injective chain map
from $C_i(S)$ to $C_i(S')$. We shall identify $C_{\bullet}(S)$ as a sub-chain complex of
$C_{\bullet}(S')$ via the function $j$.
Also, note that $z \vee x$ exists, so $C_x(z)$ is a cell of $S'$.
Thus, the $r$--cells of $S'$ are the $r$--cells of $S$ and the cones on the $(r-1)$--cells of $\tilde{S}$. (Recall that $\tilde{S}(r) = S(r)$ for $r \geq 0$ and
$\tilde{S}(-1) = \lbrace \emptyset \rbrace$.)
\par
As $\Delta(z) \cap U(x) = \emptyset$ for each facet $z$ of $y$, using equation \eqref{scone},
the boundary of a cone is
given by
\begin{equation*}
\partial([C_x(z)]) = [z] - \sum_{w \in \Delta z } s(z,w) [C_x(w)].
\end{equation*}
Let $C_{\bullet}(\tilde{S})$ be the chain complex $C_{\bullet}(S)$ augmented
by ${C}_{-1}(S) = \ZZ[\emptyset]$:
\begin{equation*}
\dotsb  \to C_i(\tilde{S}) \to C_{i-1}(\tilde{S}) \to \dotsb \to
C_1(\tilde{S}) \to C_0(\tilde{S}) \to C_{-1}(\tilde{S}) \to 0,
\end{equation*}
where the boundary map $C_0(\tilde{S}) \to C_{-1}(\tilde{S})$
sends $[x]$ to $[\emptyset]$ for each vertex $x$ of $S$.
The $i$--th homology of this chain complex will be denoted by $H_i(\tilde{S})$ for $i \geq -1$.
Let 
\begin{equation*}
h_i \co C_i(\tilde{S}) \to C_{i+1}(S')
\end{equation*}
 be the linear map induced by
$[z] \mapsto [C_x(z)]$. From the above formula for the boundary of a cone, one gets
$(h \circ \partial + \partial \circ h )([z]) = [z]$, which implies part (a).
\par
Recall that, we have identified $C_i(S)$ as a sub-complex of $C_i(S')$, via the inclusion $j$.
The function $h_i$ above induces a map $\bar{h}_i \co C_i(\tilde{S}) \to C_{i+1}(S')/C_{i+1}(S)$,
satisfying $\bar{h} \circ \partial + \partial \circ \bar{h} = 0$, showing that
$(-1)^i \bar{h}_i \co C_i(\tilde{S}) \to C_{i+1}(S')/C_{i+1}(S)$ is a chain map. The map
$\bar{h}_i$ is a bijection on the level of chains, since
$C_{i+1}(S') = C_{i+1}(S) \oplus h_i(C_i(\tilde {S}))$ as abelian groups.
So the chain complex $C_{\bullet}(S')/C_{\bullet}(S)$ is isomorphic to $C_{\bullet -1}(\tilde{S})$.
One has the following exact sequence of chain complexes:
\begin{equation*}
0 \rightarrow C_i(S) \rightarrow C_i(S') \rightarrow C_i(S')/C_i(S) \simeq
C_{i-1}(\tilde{S}) \rightarrow 0.
\end{equation*}
By taking the long exact sequence of homology groups, one gets
$H_i(S') = 0$ for $i \geq 2$, since $H_i(S) = 0$ and $H_{i-1}(\tilde{S}) = H_{i-1}(S) = 0$.
The end of this long exact sequence has the form,
\begin{equation*}
0 \rightarrow H_1(S) \to H_1(S') \to H_0(\tilde{S}) \to H_0(S) \to H_0(S') \to H_{-1}(\tilde{S}) \to 0.
\end{equation*}
By remark \ref{orientablecellremark}(3), $H_0(S) \simeq \ZZ$ is generated by the class of
any vertex of $S$. So
\begin{equation*}
\partial(C_1(S)) = \op{span} \lbrace [u] - [v]: u, v \in S(0)\rbrace.
\end{equation*}
So $\partial(C_1(S))$ is the kernel of the map $C_0(S) \to C_{-1}(S)$. Thus
$H_0(\tilde{S}) = 0$. Also $H_{-1}(\tilde{S}) = 0$. It follows that
$H_1(S') \simeq H_1(S) = 0$ and $H_0(S') \simeq  H_0(S) \simeq \ZZ$.
This finishes the proof of part (b). Part (c) follows from part (b).
\end{proof}
%
%
%
\begin{lemma}
(a) Let $S$ be a c.c.c with orientable cells and $x \in S$. Assume that $S$
is a star around $x$, that is, $\st_S(x)  = S$. Then $S^x$ is acyclic.
\newline
(b) Let $X$ be a c.c.c with orientable cells and $x \in X$.
Then $ \st(C_x(\emptyset)) \subseteq X^{x}$ is acyclic.
\label{acyclicSx}
\end{lemma}
As remarked in \ref{stellerremark}(1), there is a canonical isomorphism, $\st_{X^x}(C_x(\emptyset)) \simeq \st_X(x)^x$.
So part (b) follows from part (a). The proof of part (a), given in appendix \ref{acyclicSxpf}, is similar to the proof
of lemma \ref{cellstayacyclic}.
%
%
\begin{lemma}
Let $S$ be a c.c.c with orientable acyclic cells.
If $S$ is a star around $x$, then $S$ is acyclic.
In particular $\st(x)$ is acyclic for all $x \in S$.
(For the proof, it is important to note that we do not assume $S$ to be
equidimensional or nonsingular).
\label{acyclicS}
\end{lemma}
\begin{proof}
Let $\dim(S) = n$.
 If $x$ is a maximal cell of $S$, then  $S = \st(x) = \cl(x)$ is acyclic, by assumption.
For a non-maximal cell $x$, let $t_1, \dotsb, t_k$ be the maximal cells above $x$
arranged in decreasing order of rank, that is, $\rk(t_1) \geq \rk(t_2) \geq \dotsb \geq \rk(t_k)$.
Let $\rho_S(x) = \rk(t_1) - \rk(x)$. The proof is by induction on $\rho_S(x)$.
\par
Though logically it is not necessary, we first prove the lemma in the case
$\rho_S(x) = 1$, to illustrate the idea. Since $x$ is not a maximal cell,
one has $\rk(t_i) = \rk(x) +1$. In other words, $\nabla x = \lbrace t_1, \dotsb, t_k \rbrace $.
By induction on $j$, we show that $T_j = \cl \lbrace  t_1, \dotsb, t_j \rbrace$
is acyclic. The case $j = 1$ is a part of assumption. Assume now, that $T_{j-1}$ is acyclic.
Since $T_j = T_{j-1} \cup \cl(t_j)$ and $T_{j-1} \cap \cl(t_j) = \cl(x)$, one has the
following exact sequence of chain complexes:
\begin{equation*}
0 \to C_{\bullet}(\cl(x)) \xrightarrow[]{p} C_{\bullet}(T_{j-1}) \oplus C_{\bullet}(\cl(t_j))
\xrightarrow[]{q} C_{\bullet}(T_j) \to 0,
\end{equation*}
where $p(\lambda) = (\lambda, -\lambda)$ and $q(\mu, \sigma) = \mu + \sigma$.
By taking the long exact sequence of homology groups, one gets $H_i(T_j) = 0$ for $i \geq 2$.
Further, looking at the end of the long exact sequence, one has,
\begin{equation*}
0 \to H_1(T_j) \to H_0(\cl(x)) \xrightarrow[]{H_0(p)} H_0( T_{j-1}) \oplus H_0( \cl(t_j)) \to H_0(T_j) \to 0.
\end{equation*}
Let $v$ be any vertex of $x$. Then, by remark \ref{orientablecellremark}(3),
$[v]$ generates $H_0(\cl(x))$ and $H_0(\cl(t_j))$. The map 
$H_0(p) \co H_0(\cl(x)) \to H_0(T_{j-1}) \oplus H_0(\cl(t_j))$ sends $[v]$ to $([v],-[v])$.
Since $-[v] \in H_0(\cl(t_j))$ is non-zero, the map $H_0(p)$ is injective. It follows
that $H_1(T_j) = 0$ and $H_0(T_j) \simeq \ZZ$.
This completes the proof for $\rho_S(x) = 1$.
\par
Now, let $\rho_S(x) = r$. Assume that the lemma is true for $\rho_S(x) < r$.
By induction on $j$, we show that $T_j = \cl \lbrace  t_1, \dotsb, t_j \rbrace$
is acyclic. The case $j = 1$ is again a part of assumption. Now assume that $T_{j-1}$ is acyclic.
One has $T_j = T_{j-1} \cup \cl(t_j)$. Let $K = T_{j-1} \cap \cl(t_j)$
\footnote{Observe that $K$ is a c.c.c with orientable acyclic cells, but $K$ need not be
non-singular or equidimensional.}.
The c.c.c $K$ is a star around $x$ with $\op{dim}(K)< \rk(t_j)$, so
\begin{equation*}
\rho_K(x) < \rk(t_j) - \rk(x) \leq  r.
\end{equation*}
Since the lemma is assumed to be true for $\rho_S(x) < r$, we get that $K$ is acyclic.
As before, one has the exact sequence
\begin{equation*}
0 \to C_{\bullet}(K) \xrightarrow[]{p} C_{\bullet}(T_{j-1}) \oplus C_{\bullet}(\cl(t_j))
\xrightarrow[]{q} C_{\bullet}(T_j) \to 0.
\end{equation*}
The result follows by taking the long exact sequence of homology groups.
\end{proof}
\begin{proposition}
Assume that $X$ is a c.c.c with orientable and acyclic cells. Let $x \in X$. Then the
map $H_{\bullet}(\varphi) \co H_{\bullet}(X) \to H_{\bullet}(X^x)$, defined in \ref{StoSx},
is an isomorphism.
\label{inv}
\end{proposition}
\begin{proof}
From lemma \ref{StoSx} we have a chain map $\varphi \co C_{\bullet}(X) \to C_{\bullet}(X^x)$.
Let $S = \st_X(x)$. We shall identify $S^x$ as a sub-complex of $X^x$ via the identification
$S^x \simeq \st_{X^x}(C_x(\phi))$ given in \ref{stellerremark}(1).
The map $\varphi$ fits into the following commutative diagram of chain complexes:
\newline
\centerline{
\xymatrix{
0 \ar@{->}[r] & C_i(S) \ar@{->}[r] \ar@{->}[d]_{\varphi\vert_S} & C_i(X) \ar@{->}[d]^{\varphi} \ar@{->}[r] 
& C_i(X)/C_i(S) \ar@{->}[d]^{\wr}_{\bar{\varphi}} \ar@{->}[r] & 0  \\
0 \ar@{->}[r] & C_i(S^x) \ar@{->}[r] & C_i(X^x) \ar@{->}[r]  & C_i(X^x)/C_i(S^x) \ar@{->}[r] & 0
}}
The horizontal maps on the right are the quotient maps. 
One checks from the definitions that both $C_i(X)/C_i(S)$ and $C_i(X^x)/C_i(S^x)$ 
can be identified with the free abelian group
on the cells of $(X \setminus S)$ and the map $\bar{\varphi}$ acts as identity on these cells.
Thus $\bar{\varphi}$ is a chain isomorphism, so
$H_{\bullet}(\bar{\varphi}) \co H_{\bullet}(X,S) \to H_{\bullet}(X^x,S^x)$
is an isomorphism.
\par
Next, note that $S$ and $S^x$ are acyclic by lemma \ref{acyclicS} and \ref{acyclicSx} respectively.\footnote{
We can conclude that $S^x$ is acyclic without using lemma \ref{acyclicSx} as follows. By 
lemma \ref{cellstayorientable} and \ref{cellstayacyclic}, the cells of $X^x$ are
orientable and acyclic. So lemma \ref{acyclicS} implies $\st_{X^x}(C_x(\emptyset))$ is
acyclic. But $\st_{X^x}(C_x(\emptyset)) \simeq S^x$.}
It follows that $H_{\bullet}(\varphi \vert_S)$ is an isomorphism. Taking the diagram of homology
groups corresponding to the above commutative diagram of chain complexes and applying the five lemma,
it follows that $H_{\bullet}(\varphi) \co H_{\bullet}(X) \to H_{\bullet}(X^x)$ is an isomorphism. 
\end{proof}
%
%
%
%
\section{Barycentric subdivision of a c.c.c}
%
%
%
%
\label{secbarycentric}
Recall, from \ref{orientationdef}, the definition of the barycentric subdivision of a c.c.c $S$,
denoted by $S^{(1)}$. 
\begin{remark}
 If $S$ is equidimensional, of dimension $n$, then the same holds for $S^{(1)}$.
The $n$--cells of $S^{(1)}$ correspond to the flags in $S$. The other cells of $S^{(1)}$ correspond to
totally ordered subsets of $S$, that is, ``partial flags" in $S$.
If $S$ is non-singular, then it is easy to see that $S^{(1)}$ is non-singular.
\end{remark}
%
%
\begin{lemma}
Each cell of $S^{(1)}$ is flag connected and has an orientation such that,
for $\gamma = \lbrace x_0 > x_1 > \dotsb > x_r \rbrace \in S^{(1)}(r)$,
one has $s(\gamma, \gamma \setminus \lbrace x_i \rbrace) = (-1)^i$.
\label{S1cellorient}
\end{lemma}
\begin{proof}
The lemma follows from \ref{simphom}, once we note that
there is a compatible family of total ordering on the vertices of each cell
$\gamma \in S^{(1)}$, coming from the partial order on $S$.
\end{proof}
%
%
%
\begin{lemma}
Let $S$ be a c.c.c with orientable cells.
For each cell $x \in S$, choose an orientation $\omega_x \co \cF(x) \to \lbrace \pm 1 \rbrace $.
Choose orientations on the cells of $S^{(1)}$ as prescribed by lemma \ref{S1cellorient}.
If $x \in S(r)$, then a flag $\gamma \in \cF(x )$
determines an $r$--cell in $S^{(1)}$ and thus an $r$--chain $[\gamma]$.
There is a chain map $\Phi \co C_{\bullet}(S) \to C_{\bullet}(S^{(1)}) $, given by
\begin{equation*}
\Phi([x]) = \sum_{\gamma \in \cF(x)} \omega_x(\gamma) [\gamma].
\end{equation*}
\label{StoS1}
\end{lemma}
\begin{proof}
To check that $\Phi$ is a chain map, we first calculate $\partial_{S^{(1)}} ( \Phi[x] )$.
\begin{equation*}
\partial_{S^{(1)}} ( \Phi [x] )
=  \sum_{\gamma \in \cF(x)} \omega_x(\gamma) \partial_{S^{(1)}}  [\gamma]
=  \sum_{\gamma \in \cF(x)}  \omega_x(\gamma) \sum_{ \xi \in \Delta \gamma } s( \gamma, \xi) [\xi].
\end{equation*}
Consider a ``partial flag" $\xi$ appearing in the final expression. Suppose $\xi$ is of the form
$\lbrace x = x_0 > x_1 > \dotsb > x_{i-1} > x_{i+1} > \dotsb > x_r \rbrace$ for some $i >0$, where
$x_j\in S(r - j)$. Then there are two adjacent flags $ \gamma_+$ and  $\gamma_-$ in $\cF(x)$,
such that $\xi$ is a face of $\gamma_{\pm}$. We have $\omega_x( \gamma_+) = - \omega_x(\gamma_-)$ and
$s( \gamma_+, \xi) = s( \gamma_-, \xi) = (-1)^i$ (by lemma \ref{S1cellorient}).
So, in the expression for $\partial_{S^{(1)}} ( \Phi [x] )$, the coefficient of $[\xi]$ vanishes.
\par
Let $\xi$ be a ``partial flag" that is not of the above form.
Then $\xi$ is of the form
$\lbrace x_1 > x_2 > \dotsb > x_r \rbrace$, where $x_j$ is a cell below $x$ of rank  $(r - j)$.
That is, $\xi$ is a flag in $\cl( \Delta x)$.
The only flag $\gamma \in \cF(x)$, that has $\xi$ as a face,
is $ \gamma = \lbrace x = x_0  > x_1 > \dotsb > x_r \rbrace$.
Lemma \ref{S1cellorient} implies $s( \gamma, \xi) = 1$. It follows that
\begin{equation*}
 \partial_{S^{(1)}} ( \Phi [x] )
 =  \sum_{\xi \in \cF(\cl( \Delta x)) }  \omega_x(\xi \cup \lbrace x \rbrace ) [\xi]
 =  \sum_{y \in \Delta x } s(x,y) \sum_{\xi \in \cF( y) } \omega_y(\xi ) [\xi]
 =  \Phi( \partial_S[x]).
\end{equation*}
\par
So $\Phi$ induces a map $H_i(\Phi)\co H_i(S) \to H_i(S^{(1)})$.
\end{proof}
%
%
Suppose $S$ is a c.c.c of dimension $n$.
Let $y_1, \dotsb, y_N$ be an ordering of all the cells of $S$ of rank at-least one, such that
$\rk(y_1) \geq \rk(y_2) \geq \dotsb \geq \rk(y_N)$. 
We shall now prove that the first barycentric subdivision of $S$ can be obtained
by taking successive stellar subdivision at $y_1, y_2, \dotsb, y_N$, in that order.
Because of lemma \ref{commute}, it does not matter how the cells having the same rank are ordered.
(See proposition 2.23 of \cite{DK:CA} for the same result for simplicial complexes.)
We shall use the following abbreviation and convention:
\begin{equation*}
C_{w_j \dotsb w_1}(v) = C_{w_j} (C_{w_{j-1}}( \dotsb C_{w_1}(v) )). \text{\; If $j=0$, then \;}
C_{\emptyset}(v) = v.
\end{equation*}
%
%
\begin{lemma}
\label{baryviasteller}
 Let $S$ be a manifold--like c.c.c of dimension $n$ with orientable cells.
Let $x^r_1, x^r_2, \dotsb, x^r_{k_r}$ be the set of $r$--cells of $S$.
Starting with $T_{n+1} = S$, we shall define $T_r$ for $n+1 \geq r \geq 1$,
by backward induction on $r$. Having defined $T_{n+1}, T_n, \dotsb, T_{r+1},$ we
claim that each $r$--cell of $S$ survive as a cell of $T_{r+1}$
and we define
\begin{equation*}
T_r = T_{r+1}^{x_1^r  x_2^r \dotsb x_{k_r}^r }.
\end{equation*}
Then one has the following: 
\newline
$(\mathcal{A}(r))$ The cells of $T_r$ have the form $ C_{u_j u_{j-1} \dotsb u_1}( v )$,
where $0  \leq j \leq n-r +1$, $u_i \in S$ and $v \in S \cup \lbrace \emptyset \rbrace$.
More precisely,
\begin{equation*}
 T_r = \bigcup_{j = 0}^{n-r+1} \lbrace C_{u_j u_{j-1} \dotsb u_1}( v )  \colon
v < u_j < u_{j-1} \dotsb < u_1,  \rk(v) < r, \rk(u_{j - i}) \geq r +  i
\rbrace.
\end{equation*}
$(\mathcal{B}(r))$ The cells greater than or equal to $C_{w_j \dotsb w_1}(t)$ in $T_{r+1}$ are the cells of the form
$C_{y_k \dotsb y_1}(v)$ where $t \leq v$ and $\lbrace w_j < w_{j-1} < \dotsb < w_1 \rbrace$
is an ordered subset of $\lbrace y_{k} < y_{k-1} < \dotsb < y_1 \rbrace$.
\newline
$(\mathcal{C}(r))$ Consider $t \in S(r)$ as a cell of $T_{r+1}$. The cells of $T_{r+1}$ that are greater than or equal
to $t$ are those of the form $C_{u_j u_{j-1} \dotsb u_1} (t) $, $ j \geq 0$. Thus, if $t$ and  $t'$
are two distinct $r$--cells of $S$, then $U_{T_{r+1}}(t) \cap U_{T_{r+1}}(t') = \emptyset$.
Consequently $(T_{r+1})^{t  t'} \simeq (T_{r+1})^{t' t}$.
\newline
$(\mathcal{D})$ The c.c.c $T_1$ is canonically isomorphic to the first barycentric subdivision $S^{(1)}$. Under this isomorphism,
The cell $C_{v_1 v_2 \dotsb v_r}(v_0) \in T_1$ corresponds to the cell
$\lbrace  v_0 < v_1 < \dotsb < v_r \rbrace \in S^{(1)}$. 
\end{lemma}
In the statement $\mathcal{A}(r)$ 
The proof is given in \ref{baryviastellerpf}. However, it is best to work out a
few examples in dimension 2 and 3 to convince oneself of the validity of the statement.
%
%
\begin{proposition}
Let $S$ be a manifold--like c.c.c of dimension $n$, with orientable and acyclic cells. Then
$H_{\bullet}(S) \simeq H_{\bullet}(S^{(1)})$.
\label{baryinv}
\end{proposition}
\begin{proof}
By lemma \ref{baryviasteller}, the first barycentric subdivision $S^{(1)}$ is obtained from $S$
by a sequence of successive stellar subdivisions.
The property of having orientable and acyclic cells, is preserved under stellar subdivision,
by lemma \ref{cellstayorientable} and \ref{cellstayacyclic} respectively.
The result now follows from repeated application of proposition \ref{inv},
which says that, for c.c.c's with acyclic orientable cells, homology is invariant under stellar
subdivision.
\end{proof}
%
%
\begin{remark}
\label{baryinvf}
We can refine proposition \ref{baryinv}, as follows.
Let $y_1, \dotsb, y_N$ be a list of all the cells of $S$ in decreasing order of rank. 
Let $\varphi^{\circ}$ be the composite of the chain maps given below:
\begin{equation*}
C_{\bullet}(S) \to C_{\bullet}(S^{y_1}) \to  C_{\bullet}(S^{y_1 y_2}) \to \dotsb \to 
C_{\bullet}( S^{y_1 y_2 \dotsb y_N}) \simeq  C_{\bullet}( S^{(1)})
\end{equation*}
where all but the last chain map is obtained from lemma \ref{StoSx}
and the last isomorphism is a consequence of lemma \ref{baryviasteller}.
It follows from lemma \ref{inv}, that $\varphi^{\circ}\co C_{\bullet}(S) \to C_{\bullet}(S^{(1)})$
induces isomorphisms of homology groups.
On the other hand, lemma \ref{StoS1} gives us another chain map 
$\Phi \co C_{\bullet}(S) \to C_{\bullet}(S^{(1)})$.
One can check that 
\begin{equation}
\Phi_j = \pm \varphi^{\circ}_j\co C_j(S) \to C_j(S^{(1)}).
\label{Phiphi}
\end{equation} 
(A proof of equation \eqref{Phiphi} is given in appendix \ref{Phiphipf}).
From equation \eqref{Phiphi} it follows that
$H_{\bullet}(\Phi) \co H_{\bullet}(S) \to H_{\bullet}(S^{(1)})$ is an isomorphism.
\par
There is a somewhat confusing issue here, that needs an explanation.
It follows from \ref{StoSx} and \ref{StoS1} that both $\varphi^{\circ}$
and $\Phi$ commute with the boundary maps. However, the maps
$\varphi^{\circ}$ and $\Phi$ only agree up-to sign. 
The solution to this apparent contradiction is the following observation.
To show that $\Phi$ (resp. $\varphi^{\circ}$) is a chain map we must orient the cells of 
$S^{(1)}$ as prescribed by lemma \ref{S1cellorient} (resp. repeated use of lemma
\ref{cellstayorientable}). These two sets of orientations on the cells of
$S^{(1)}$ do not agree. So the two boundary maps on $S^{(1)}$, with respect
to which $\varphi^{\circ}$ and $\Phi$ are shown to be chain maps, are different.
\end{remark}
\section{Poincare duality}
\label{secthm}
\begin{lemma}
Let $S$ be an orientable, manifold--like c.c.c of dimension $n$. Assume that each cell of $S$ and $S^{\circ}$
is flag--connected. Then \\
(a) $S^{(1)} = (S^{\circ})^{(1)}$.\\
(b) $H_i( S) \simeq H^{n-i}(S^{\circ})$.
\label{SSop}
\end{lemma}
\begin{proof}
Proof of part (a) is clear from the definitions.
\par
Proof of part (b) is like the classical proof of Poincare duality theorem,
by relating homology and cohomology using dual cell decompositions (for example, see \cite{GH:AG}, pages 53--55).
Since $S$ is orientable, manifold--like, of dimension $n$, so is $S^{\circ}$
(by \ref{Sop}).
Since $S$ is orientable and each cell of $S$ is flag--connected,
the first remark in \ref{orientablecellremark} implies that each cell of $S$
is orientable. The same remark holds for $S^{\circ}$.
\par
Recall that the flags in $\cF(x) = \cF(\cl_S(x))$ are called the flags below $x$
and the flags in $\cF(x^{\circ}) = \cF(\cl_{S^{\circ}}(x^{\circ}))$ are called the flags above $x$. 
Suppose $y = x$ or $y \in \Delta x$ and we are given a flag $\gamma_2$  above $x$
and a flag $\gamma_1$ below $y$. 
Then, putting together $\gamma_1$ and $\gamma_2$, with the partial order on $\gamma_2$ reversed,
one obtains a flag in $S$, which we shall denote by
$\gamma_1 \cup \gamma_2^{\circ} $.
\par
Let $\omega$ be an orientation on $S$ and $\omega^{\circ}$ be the corresponding
orientation on $S^{\circ}$.
For each $x \in S$, choose an orientation $\omega_x$ on $\cl_S(x)$
such that, if $x$ is a maximal cell, then $\omega_x$ is the restriction
of $\omega$ to $\cl_S(x)$.
Define an orientation $\omega_x^{\circ}$ on $\cl_{S^{\circ}}(x^{\circ})$ as follows.
Given a flag $\gamma_2 \in \cF(x^{\circ})$, choose flag $\gamma_1$ below $x$ and
define
\begin{equation*}
 \omega_x^{\circ}(\gamma_2) = \omega(\gamma_1 \cup \gamma_2^{\circ} )/\omega_x(\gamma_1).
\end{equation*}
The definition of $\omega_x^{\circ}$ does not depend on the choice of $\gamma_1$, because
the adjacency graph of flags below $x$, is connected. Further, if $\gamma_2$ and $\tilde{\gamma}_2$
are adjacent flags above $x$, then $\gamma_1 \cup \gamma_2^{\circ}$ and $\gamma_1 \cup \tilde{\gamma}_2^{\circ}$
are adjacent flags in $S$. It follows that
\begin{equation*}
\omega_x^{\circ}(\gamma_2) = \omega( \gamma_1 \cup \gamma_2^{\circ})/\omega_x(\gamma_1)
= - \omega(\gamma_1 \cup \tilde{\gamma}_2^{\circ}) /\omega_x(\gamma_1) = -\omega_x^{\circ}(\tilde{\gamma}_2),
\end{equation*}
showing that $\omega_x^{\circ}$ is an orientation on $\cl(x^{\circ})$.
\par
Now suppose that $y$ is a face of $x \in S$. Pick a flag $\gamma_1$ below $y$ and a flag $\gamma_2$ above $x$,
and let $\gamma = \gamma_1 \cup \gamma_2^{\circ}$ be the flag in $S$, obtained by putting them together
\footnote{If $\gamma_2 = \lbrace x_{n-r-1}^{\circ} <^{\circ} \dotsb <^{\circ} x_0^{\circ} = x \rbrace$ 
and $\gamma_1 = \lbrace y = y_0 > \dotsb > y_r  \rbrace$ then
putting them together one gets the flag
$\gamma = \lbrace x_{n-r-1} > x_{n-r-2} > \dotsb > x_1  > x_0 > y_0 > y_1 > \dotsb > y_r \rbrace$.}.
Then one has
\begin{equation}
 s^{\circ}(y^{\circ}, x^{\circ})
= \frac{\omega_y^{\circ}( \gamma_2 \cup \lbrace y^{\circ} \rbrace )}{\omega_x^{\circ}(\gamma_2)}
=\frac{\omega(\gamma)/\omega_y(\gamma_1)}{\omega(\gamma)/\omega_x(\gamma_1 \cup \lbrace x \rbrace )}
= s(x,y).
\label{sso}
\end{equation}
Consider the map 
$* \co C_i(S) \to C_{n-i}(S^{\circ})$ given by $*[x] = [x^{\circ}]$.
The equation \eqref{sso} shows that
\begin{equation*}
*(\partial[\sigma] ) = \delta ( *[\sigma] ), \text{\; for \;} \sigma \in C_i(S).
\end{equation*}
So the map $*$ is an isomorphism between the chain complexes
$(C_i(S), \partial)$ and $(C_{n-i}(S^{\circ}), \delta)$.
\end{proof}
\begin{theorem}
Suppose $S$ is an orientable, manifold--like c.c.c of dimension $n$. Assume that each cell of $S$
and $S^{\circ}$ is flag--connected and acyclic. Then $H_i( S) \simeq H^{n-i}(S)$.
\label{duality}
\end{theorem}
\begin{proof}
As $S$ is $n$ dimensional, manifold--like and orientable, the same holds for $S^{\circ}$.
Since $S$ is orientable and each cell of $S$ is flag--connected,
the first remark in \ref{orientablecellremark} implies that each cell of $S$ is orientable.
The same remark holds for $S^{\circ}$. So each cell of $S$ and $S^{\circ}$ is orientable and acyclic.
By proposition \ref{baryinv} the homology of $S$ and $S^{\circ}$ are invariant under
barycentric subdivision. But the barycentric subdivision of $ S$ and $S^{\circ}$ are identical
(see lemma \ref{SSop} (a)). It follows that
\begin{equation*}
H_i(S) \simeq H_i(S^{(1)} ) \simeq H_i(( S^{\circ})^{(1)} ) \simeq H_i(S^{\circ}).
\end{equation*}
Since $(S^{\circ})^{\circ} = S$, By lemma \ref{SSop}(b), we have $H_i(S^{\circ} ) \simeq H^{n-i}(S)$.
\end{proof}
%
%
\section{Miscellaneous remarks}
%
%

\begin{topic}{\bf Intersection pairing and integration:}
Let $S$ be an orientable, manifold--like c.c.c of dimension $n$.
Note that one has a tautological pairing,
\begin{equation*}
C_i(S) \times C_{n-i} (S^{\circ}) \to \ZZ,
\end{equation*}
obtained by linearly extending $\langle [x], [z^{\circ}] \rangle = \chi(x = z)$,
where $\chi(\cdot)$ is the indicator function.
Let $x \in S(i+1)$ and $z \in S(i)$. Using equation \eqref{sso}, one has,
\begin{equation*}
\langle \partial [x] , [z^{\circ}] \rangle
= \sum_{ y \in \Delta x} s(x, y) \chi(y = z)
= s(x, z) \chi(z \in \Delta x)
= s^{\circ}(z^{\circ}, x^{\circ}) \chi( x^{\circ} \in \nabla z^{\circ} )
=\langle [x] , \partial [z^{\circ}] \rangle.
\end{equation*}
By linearly extending, one gets,
\begin{equation}
\langle \partial \sigma , \tau \rangle = \langle \sigma , \partial \tau \rangle,
\label{sa}
\end{equation}
for $\sigma \in C_{i+1}(S)$ and $ \tau \in C_{n-i}(S^{\circ})$.
The pairing between chains and co-chains restricts to give a pairing between $i$--cycles of $S$
and $(n-i)$--cycles of $S^{\circ}$. Equation \eqref{sa}
shows that the pairing between cycles descends to a pairing between the homology groups,
\begin{equation*}
H_i(S) \times H_{n-i}(S^{\circ}) \to \ZZ, \text{\;\; denoted by \;\;}
(\sigma , \tau) \mapsto \sigma \pitchfork \tau.
\end{equation*}
This is the {\it intersection pairing}.
From lemma \ref{SSop}, we have an isomorphism $* \co H_i(S) \to H^{n-i}(S^{\circ})$.
Let us also denote the inverse isomorphism by $*$.
Using the duality $*$ and the intersection pairing, we get the {\it integration pairing}:
\begin{equation*}
 \int\co H_i(S) \times H^i(S) \to \ZZ, \text{\;\;defined by\;\;}
\int_{\sigma} \omega = \sigma \pitchfork * \omega.
\end{equation*}
An immediate consequence of equation \eqref{sa} is Stoke's theorem:
$\int_{\partial \sigma} \omega = \int_{\sigma} \delta \omega$.
\end{topic}
%
%
%
\begin{topic}
{\bf Functoriality of homology groups: }
\label{functor}
Let $\msf{Cat}$ be the category of small categories and
let $N$ be the nerve functor defined from $\msf{Cat}$ to the category of
simplicial sets. 
Let $\msf{CCC}$ be the category whose objects
are combinatorial cell complexes and the morphisms are order preserving maps of underlying posets,
or in other words, continuous maps of the underlying finite topological spaces.
Considering a partially ordered set as a category with only one morphism between any two objects,
we can view $\msf{CCC}$ as a full subcategory of $\msf{Cat}$.
Thus, given a c.c.c $X$, we get a simplicial set $N(X)$,
whose $r$-simplices are 
\begin{equation*}
N(X)_r = \lbrace (x_0, x_1, \dotsb, x_r) \colon x_0 \geq x_1 \geq \dotsb \geq x_r , x_j \in X\rbrace
\end{equation*}
and the $j$--th face map is given by $\partial_j(x_0, \dotsb, x_r) = (x_0, \dotsb, x_{j-1}, x_{j+1}, \dotsb, x_r)$.
\par
Let us recall the definition of the normalized homology groups of the simplicial set $N(X)$.
The boundary map $\partial\co \ZZ[N(X)_r] \to \ZZ[N(X)_{r-1}]$ is
obtained by linearly extending $\partial x = \sum_j (-1)^j \partial_j x$. The homology of
the simplicial set $N(X)$ is the homology of the chain complex $(\ZZ[N(x)_{\bullet}], \partial)$.
The chains supported on degenerate cells,\footnote{$(x_0 \geq x_1 \geq \dotsb \geq x_r)$ is a
degenerate cell of $N(X)$ if $x_j = x_{j+1}$ for some $j$.}
form a sub-complex of the above chain complex and the homology groups of the quotient chain complex
are the normalized homology groups of $N(X)$.
It is classically known\footnote{See 10.6 of \cite{EZ:SS}. Simplicial sets were first defined
in this article under the name ``complete semi simplicial complexes".} that the quotient maps
on chains induce canonical isomorphisms from the homology groups of a simplicial set to the
normalized homology groups. 
\par
Let $\gamma = \lbrace x_0 > \dotsb > x_r \rbrace$ be an $r$--cell of $X^{(1)}$.
From lemma \ref{S1cellorient}, recall that the boundary map for the chain complex of the c.c.c $X^{(1)}$, is given by
\begin{equation*}
\partial [\gamma] = \sum_j s(\gamma, \gamma \setminus \lbrace x_j \rbrace) [\gamma \setminus \lbrace x_j \rbrace]
=\sum_j (-1)^j [x_0 > \dotsb > x_{j-1} > x_{j+1} >  \dotsb > x_r].
\end{equation*}
So the inclusion $X^{(1)} \hookrightarrow N(X)$, induces a chain map
from $(C_{\bullet}(X^{(1)}), \partial) \to (\ZZ[N(X)_{\bullet}], \partial)$,
which, after quotienting out on the right by the group generated by the
degenerate cells, becomes an isomorphism, since
the $r$--cells of $X^{(1)}$ are precisely the non-degenerate $r$--cells of $N(X)$.
It follows that the homology of the c.c.c $X^{(1)}$ is canonically isomorphic
to the normalized homology of the simplicial set $N(X)$,
which is canonically isomorphic to the homology of $N(X)$.
\par
Let $X$ and $Y$ be combinatorial cell complexes with orientable cells.
Given a continuous map $f\co X \to Y$ of finite spaces, it is not in general clear
how to get a map between the cellular homology groups that we defined in section \ref{sechomcohom}.
However, consider the subcategory $\msf{CCC}_a \subseteq \msf{CCC}$, consisting of manifold--like
combinatorial cell complexes with orientable and acyclic cells.
Let $X$ be an object of $\msf{CCC}_a$. From \ref{baryinvf}, one has a canonical isomorphism
$H_{\bullet}(\Phi)\co H_{\bullet}(X) \to H_{\bullet}(X^{(1)})$.
Composing with the canonical isomorphism $H_{\bullet}(X^{(1)}) \to H_{\bullet}(N(X))$,
one gets a canonical isomorphism $ \Phi^X_*\co H_{\bullet}(X) \to H_{\bullet}(N(X))$, for
each object $X$ of $\msf{CCC}_a$.
Thus, given a morphism $f \co X \to Y$ in $\msf{CCC}_a$, one gets
an induced morphism of abelian groups, $H_i(f)\co H_i(X) \to H_i(Y)$, defined by
$H_i(f) = (\Phi^Y_*)^{-1} \circ H_i( N(f)) \circ \Phi^X_*$, for all $i$. 
Since $N(\cdot)$ is a functor and $H_i$ are functors on simplicial
sets, it follows that $H_i$ are functors from $\msf{CCC}_a$ to abelian groups.
\end{topic}
\begin{topic}
{\bf Infinite combinatorial cell complexes: }
In the definition of a c.c.c $(S, \leq, \rk)$, given in \ref{cccdef},
suppose we allow the poset $S$ to be infinite. The definition still makes sense.
Many of the results in this article hold for infinite $S$,
if we only assume that $\cl(x)$ is finite for all $x \in S$.
Most results hold if we assume that $S$ is finite dimensional and 
that for each $x \in S$, both $\cl(x)$ and $U(x)$ are finite.
The exact finiteness condition, that needs to be imposed on
$S$ for a particular lemma, should be clear by looking at the proof.
For the sake of clarity, we have assumed throughout that $S$ is finite. 
\end{topic}
%
%
%
%
\appendix
\section{Proofs of some lemmas}
\begin{topic}
\begin{proof}[proof of lemma \ref{stellerfirstlemma}]
\label{stellerfirstlemmapf}
(a) {\bf Axiom (1):} Recall that $y <_{S^x} z$ if and only if one of the following
three conditions hold:
$(i)$ $y \in S$, $z \in S$ and $y <_S z$, or
$(ii)$ $y \in S$, $z = C_x(z')$ and $y \leq_S z'$, or
$(iii)$ $y = C_x(y')$, $z = C_x(z')$ and $ y' <_S z'$.
In each of these cases, $\rk_{S^x}(y) < \rk_{S^x}(z)$.
\par
{\bf Axiom (2):}
Let $T$ be a subset of $S^x$ that is bounded below.
Let $\tilde{T}_N = \lbrace v \in S : C_x(v) \in T \rbrace$ and $T_{O} = T \cap S$.
If $T_O \neq \emptyset$, then any lower bound $y$ of $T$ is necessarily an old cell.
Then both $T_O$ and   $\tilde{T}_N$ are bounded below by $y$ and
$\wedge T = \wedge (T_O \cup  \tilde{T}_N )$.
On the other hand, if $T_O = \emptyset$, then $\tilde{T}_N$ is bounded below,
$C_x( \wedge(\tilde{T}_N))$ exists and is equal to $\wedge T$.
Given $ y < z$ in $S^x$, it is easy to find a cell $y' \in S^x$ such that $\rk(y')  = \rk(y) + 1$
and $y < y' \leq z$.
\par
{\bf Axiom (3):} Suppose $z \in S^{x}$ is a cell of rank at-least 1 and $u$ is an upper bound
for $\Delta z$. We need to check that $u \geq z$.
First, suppose that $z$ is an old cell. If $u$ is an old cell, then $u$ is an upper bound
for $\Delta z$ in $S$, so $u \geq z$. If $u = C_x(u')$ is a new cell, then $C_x(u') \geq y$
for all $y \in \Delta z$, which implies that $u' \geq y$ for each $y \in \Delta z$, so $ u' \geq z$
and hence, $ u = C_x(u') \geq z$.
\par
Next, suppose that $z = C_x(z')$ is a new cell. Then
\begin{equation*}
\Delta z  = \lbrace z' \rbrace \cup \lbrace C_x(v) \colon v \in \Delta z' \rbrace.
\end{equation*}
Any upper bound $u$ for $\Delta z$ must be a new cell, that is, $ u = C_x(u')$.
Now, $C_x(u') \geq z'$ implies that $ u' \geq z'$ in $S$, which in turn implies that
$u = C_x(u') \geq C_x(z') = z$.
\par
{\bf Axiom (4):} Let $y$ be a co-dimension 2 facet of $z$ in $S^x$. 
If $z$ is an old cell, then the set of cells below $z$ is the same in
$S$ and $S^x$, so there are two cells between $y$ and $z$.
If $z = C_x(z')$ and $y = C_x(y')$ are both new cells,
then the cells between $z$ and $y$ in $S^x$ are in one to one correspondence with
the cells between $z'$ and $y'$ in $S$, so there are just two of them.
Finally, suppose that $z = C_x(z')$ is a new cell and $ y$ is an old cell.
Suppose $y < w < z$. If $w$ is not a cone, then $ w = z'$. If $w = C_x(w')$ is a cone,
then $y = w'$ and hence $w = C_x(y)$. (Note that $z' \in M(x)$ and $y < z'$
implies that $y \in M(x)$, so $C_x(y)$ exists).
Hence there are two cells between $z = C_x(z')$ and $y$, namely
$C_x(y)$ and $z'$.
\newline
\par
(b) Let $S$ be equidimensional, of dimension $n$. Let $x \in S$ and $t \in M(x)$.
\newline
{\it Claim: There exists a cell $w \in M(x)$, such that $w \geq t$ and $\rk(w) = n-1$.}
\newline
{\it proof of the claim:} Let $w$ be a cell of maximal rank above $t$ in $M(x)$. Suppose, if possible,
that $\rk(w) < n-1$. If $\rk(w \vee x) > \rk(w) + 1$, then there would be
a cell strictly in between $w$ and $w \vee x$, which would contradict the maximality
of $w$. Thus $\rk(w \vee x) = \rk(w) + 1 < n$.
So there is a cell $z$, such that $ z_+ = w \vee x $ is a face of
$z$. But there is another face of $z$, call it $z_-$, between $z$ and $w$. If
$z_- \notin U(x)$, then the maximality of $w$ is contradicted. On the other
hand, if $z_- \in U(x)$, then $w = z_+ \wedge z_- \geq x$, which is again a contradiction.
This proves the claim.
\par
Let $t$ be a non-maximal cell of  $S^x$.
We need to show that there is an $n$--cell of $S^x$ above $t$. Suppose $t$ is an old cell.
If there is an $n$--cell of $S$, that is above $t$ but not above $x$, then we are done. So assume
that all the $n$--cells above $t$ are in $U(x)$. In particular $t \in M(x)$.  By the claim above, there is a $w \geq t$ in $M(x)$ of rank $n-1$. 
So $C_x(w)$ exists and is a $n$--cell in $S^x$ above $t$.
\par
Now, suppose that $t$ is a new cell, that is, $t =  C_x(t')$ for some $t' \in M(x)$.
By the claim above, there is a $w \geq t'$ such that $w \in M(x)$ and $\rk(w) = n-1$.
So $C_x(w)$ is an $n$--cell above $C_x(t')$.
\newline
\par
(c) Let $y$ be a cell of $S^x$ of rank one. If $y$ is not a cone, then
the vertices of $y$ are also not cones, so $y$ has two vertices.
Otherwise $y = C_x(y')$ for some $y' \in S(0)$. Let $z \in \Delta y$. Either $ z = C_x(\emptyset)$
or $z$ is not a cone. In the latter case $z \leq y'$ and hence $z = y'$.
\newline
\par
(d) Suppose $S$ is equidimensional, of dimension $n$.
Suppose $y$ is an old $(n-1)$--cell of $S^x$.
If $y \notin M(x)$, then the co-faces of $y$ in $S^x$ are the same as the co-faces
of $y$ in $S$, so we have nothing to prove. So assume that $y \in M(x)$.
In this situation, $C_x(y)$ is the only cone above $y$.
If $u$ is the only $n$--cell above $y$ in $S$, then one must have $u = y \vee x$, 
so $u$ is no longer a cell of $S^x$. So $C_x(y)$ is the only $n$--cell above
$y$ in $S^x$.
Now, suppose that there are two $n$--cells $u_+$ and $u_-  = y \vee x$ above $y$
in $S$. If $u_+ \in U(x)$, then one would have $y = u_+ \wedge u_- \geq x$, which
is not true. So $u_+ \notin U(x)$. So $u_+$ and $C_x(y)$ are the two $n$--cells
above $y$ in $S^x$.
\par
Now suppose $y = C_x(y')$ is a new $(n-1)$--cell of $S^x$. Let $z = C_x(z')$ be any cell above $y$.
Then $z' \vee x$ and $y' \vee x$ exists. We summarize the situation in figure \ref{butterfly}(a).
The left rhombus is in $S$ and the right rhombus is in $S^x$.\\
\begin{figure}
\begin{center}
\includegraphics{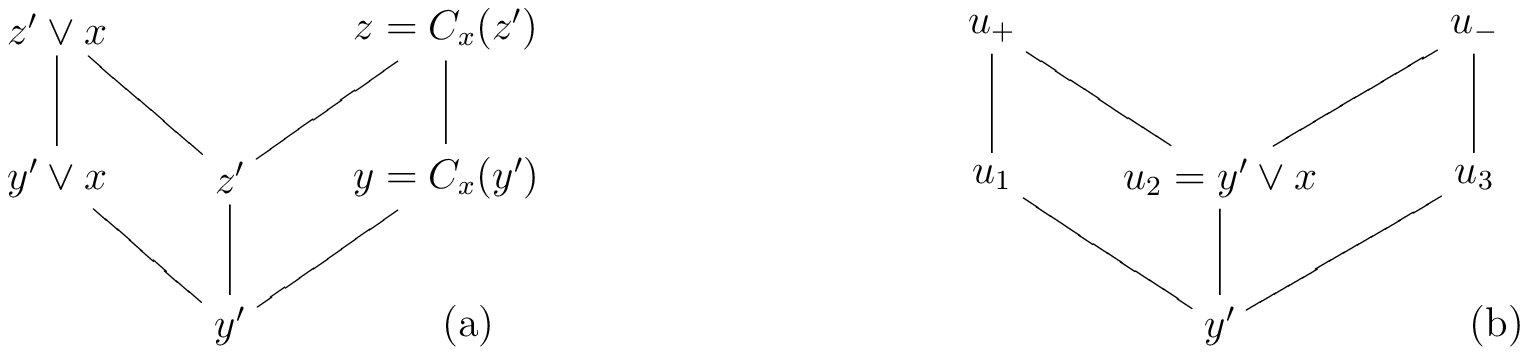}
\caption{}
\label{butterfly}
\end{center}
\end{figure}
We have to consider two cases, namely $\rk(y' \vee x) = n - 1$ and $\rk(y' \vee x) = n$.
\par
{\bf Case I : $\rk(y' \vee x) = n-1$.}
One has $C_x(z') > C_x(y')$ if and only if $z' \vee x$ is an $n$--cell above $y' \vee x$.
There are one or two $n$--cells in $S$ above $y' \vee x$.
Accordingly we have two sub-cases:
\begin{enumerate}
\item  Suppose, there is only one $n$--cell above $y' \vee x$, call it $u$. Then $ z' \vee x  = u$.
So $z'$ must be below $u$ and above $y'$. There are exactly two such cells in $S$.
One of them, namely $y' \vee x$, is not a possible choice for $z'$ since $y' \vee x \geq x$.
So there is only one choice for $z'$ and hence for $z$.
\item Suppose that there are two $n$--cells above $y' \vee x$, call them $u_+$ and $u_-$.
The purported $z'$ must be above $y'$ and below either $u_+$ or $u_-$.
By axiom (4) in the definition of a c.c.c, there are three such cells,
say $u_1, u_2$ and $u_3$, where $ u_1 < u_+$, $u_3 < u_-$ and $u_2 = u_+ \wedge u_- = y' \vee x$
(see figure \ref{butterfly}(b)).
One of them, namely $u_2$, is not a possible choice for $z'$, since $u_2 \in U(x)$.
Note that $u_1 \in U(x)$ would imply
$u_1 \wedge u_2 = y' \geq x$ which is not true. So $u_1 \notin U(x)$. For similar reason
$u_3 \notin U(x)$.
So either $z' = u_1$ implying $z = C_x(u_1)$ or $z' = u_3$ implying $z = C_x(u_3)$.
\end{enumerate}
{\bf Case II : $\rk(y' \vee x) = n$.}
In this case $z' \vee x = y' \vee x$. So the purported $z'$ must be
below $y' \vee x$ and above $y'$. There are two such cells, both in $M(x)$.
So $z'$ must equal one of them. So there are two choices for $z'$ and
correspondingly, two choices for $z$.
This finishes the proof of part (d). Part (e) now follows from (c) and (d).
\end{proof}
\end{topic}
%
%
\begin{topic}
\label{commutepf}
\begin{proof}[proof of lemma \ref{commute}]
One proceeds by induction on $k$. When $k = 1$, the lemma follows from the definition
of a stellar refinement. Assume that $X_{(k-1)}$ has the description given in the lemma.
Note that $x_k$ is an old cell of $X_{(k-1)}$. If $x_k <_{X_{(k-1)}}  C_{x_j}(v)$ for some $j < k$, then
$x_k <_X v$ and $v \vee x_j$ exists in $X$, and one has $v \vee x_j \in U_X(x_j) \cap U_X(x_k)$, which is a
contradiction. So there are no new cells of $X_{(k-1)}$ above $x_k$. If $\alpha = C_{x_j}(v)$ is a new cell in
$M_{X_{(k-1)}}(x_k)$, then $ \alpha \vee x_k$ would be a new cell of $X_{(k-1)}$ above $x_k$,
which is again impossible. So there are no new cells of $X_{(k-1)}$ in $M_{X_{(k-1)}}(x_k)$
either.
Next, observe that if $\alpha \in U_X(x_k)$ or $\alpha \in M_X(x_k)$, then $\alpha \notin \cup_{j =1}^{k-1} U_X(x_j)$,
so $\alpha$ survives as a cell in $X_{(k-1)}$. From the above discussion it follows that
\begin{equation*}
U_{X_{(k-1)}}(x_k) = U_X(x_k) \text{\;\; and \;\;}  M_{X_{(k-1)}}(x_k) = M_X(x_k).
\end{equation*}
 Hence the set
$(X_{(k-1)})^{x_k} = \lbrace C_{x_k}(v) \colon v \in \tilde{M}_X(x_k) \rbrace \cup ( X_{(k-1)} \setminus U_X(x_k))$,
matches the description of $X_{(k)}$ given in the lemma.
\par
It remains to check that the partial order on $(X_{(k-1)})^{x_k}$ matches the description given in the lemma. 
From the definition of partial order on a stellar refinement, it follows that the relation
$\alpha \leq_{X_{(k)}} \beta$ holds, if and only if one of the following three possibilities are true:
\begin{itemize}
\item Both $\alpha$ and  $\beta$ belong to  $X_{(k-1)}$ and $\alpha \leq_{X_{(k-1)}} \beta$.
By the induction hypothesis, we already know when this happens.
\item $\alpha \in X_{(k-1)} \setminus U_X(x_k)$, $ \beta = C_{x_k}(\beta')$ and
$\alpha \leq_{X_{(k-1)}} \beta'$. Here $\beta' \in M_{X_{(k-1)}}(x_k) = M_X(x_k)$
is an old cell. From the description of the partial order on $X_{(k-1)}$, it follows
that $\alpha$ must also be an old cell, so $\alpha \in X$ and $\alpha <_X \beta'$.
\item $\alpha$ and $\beta$ are of the form $\alpha = C_{x_k}(\alpha')$,
$\beta = C_{x_k}(\beta')$ for some $\alpha', \beta' \in \tilde{M}_X(x_k)$ and
$ \alpha' <_X \beta'$.
\end{itemize}
These three possibilities amount to the proposed description of the partial order
on $X_{(k)}$.
\end{proof}
\end{topic}
%
%
\begin{topic}
\label{cellstayorientablepf}
\begin{proof}[proof of lemma \ref{cellstayorientable}]
(a) One only has to show that the graph $\cF(C_x(y))$ is connected, for each $y \in M(x)$.
Let $\rk(y) = n-1$.
Let $\cF' \subseteq \cF(C_x(y))$ be the set of flags of the form
$\lbrace C_x(y_0) > C_x(y_1) > \dotsb > C_x(y_n) \rbrace$ where
$y_0 = y$ and $y_n = \emptyset$. The sub-graph of $\cF(C_x(y))$, with vertex set $\cF'$, is
isomorphic to the adjacency graph of the flags in $\cl(y)$, hence $\cF'$ is connected.
\par
Given a flag $\gamma_1$ of the form
\begin{equation*}
\gamma_1 =
\lbrace C_x(y_0) > C_x(y_1) > \dotsb  > C_x(y_i) > y_i > y_{i+2} >  \dotsb > y_n \rbrace,
\end{equation*}
one has a flag
\begin{equation*}
\gamma_2 = \lbrace C_x(y_0) > C_x(y_1) > \dotsb > C_x(y_i) >C_x(y_{i+2}) > y_{i+2} > \dotsb > y_n \rbrace
\end{equation*}
which is adjacent to $\gamma_1$ and has one more cone in it. So any flag in $\cF(C_x(y))$ is
connected to a flag consisting of all cones, that is, a flag in $\cF'$. This proves part (a).
\par
(b) Let $\gamma_1 = \lbrace a_0 > a_1 >  \dotsb > a_n \rbrace$ and
$\gamma_2 = \lbrace b_0 > b_1 > \dotsb > b_n \rbrace $ be adjacent flags in $S'$.
Assume that $a_r \neq b_r$ and $ a_j = b_j $ for $ j \neq r$.
Observe that $l(\gamma_1)$ and $l(\gamma_2)$ can differ by at-most one.
Without loss, assume that $l(\gamma_2) \geq l(\gamma_1)$.
\par
First, assume that $l(\gamma_1) = l(\gamma_2) = i$. Then the level $r$, at which $\gamma_1$
and $\gamma_2$ differs, cannot be $i$ or $(i+1)$.
It follows that
\begin{equation*}
\gamma_j = \lbrace C_x(y_0^j) >  \dotsb > C_x(y_i^j) > y_i^j > y_{i+2}^j > \dotsb > y_n^j \rbrace
\end{equation*}
for $j = 1, 2$, where $y_k^1 = y_k^2$ for all $k \neq r$ and $y_r^1 \neq y_r^2$.
So $\tilde{\gamma}_1$ and $\tilde{\gamma}_2$ are adjacent in $\cF(S)$.
It follows that
$\omega_{S'}(\gamma_1) = (-1)^i \omega_y(\tilde{\gamma}_1) = -  (-1)^i \omega_y(\tilde{\gamma}_2)
= - \omega_{S'}(\gamma_2)$.
\par
Now, assume that $l(\gamma_1) = i$ and $l(\gamma_2) = (i+1) $. The flags $\gamma_1$ and $\gamma_2$ can be
adjacent, only if they have the following form:
\begin{align*}
\gamma_1 &= \lbrace C_x(y_0) > \dotsb > C_x(y_i) &&>&  y_i         &&>&& y_{i +2} > \dotsb > y_n \rbrace \\
\gamma_2 &= \lbrace C_x(y_0) > \dotsb > C_x(y_i) &&>& C_x(y_{i+2}) &&>&& y_{i +2} > \dotsb > y_n \rbrace
\end{align*}
In this case, $\tilde{\gamma}_1 = \tilde{\gamma}_2$.
It follows that
\begin{equation*}
\omega_{S'}(\gamma_1) = (-1)^i \omega_y(\tilde{\gamma}_1) = -  (-1)^{i+1} \omega_y(\tilde{\gamma}_2)
= - \omega_{S'}(\gamma_2).
\end{equation*}
\end{proof}
\end{topic}
%
%
\begin{topic}
\label{acyclicSxpf}
\begin{proof}[proof of lemma \ref{acyclicSx}]
Let $U = U_S(x)$ and $M = M_S(x) = S \setminus U$.  Let
$\tilde{M}(i) = M(i)$ for $i \geq 0$ and $\tilde{M}(-1)= \lbrace \emptyset \rbrace$.
Let $C_{\bullet}(\tilde{M})$ be the chain complex
\begin{equation*}
0 \to C_{n-1}(\tilde{M}) \dotsb \to C_i(\tilde{M}) \to C_{i-1}(\tilde{M}) \to \dotsb \to
C_1(\tilde{M}) \to C_0(\tilde{M}) \to C_{-1}(\tilde{M}) \to 0,
\end{equation*}
where the boundary map $C_0(\tilde{M}) \to C_{-1}(\tilde{M})$ sends each vertex of $M$ to
$[\emptyset]$. Let $H_{\bullet}(\tilde{M})$ be the homology of the complex $C_{\bullet}(\tilde{M})$.
One has $H_i(\tilde{M}) = H_i(M)$ for $i \geq 1$, $H_{-1}(\tilde{M}) = 0$ and
$H_0(\tilde{M})$ is a free abelian group with $\rk(H_0(\tilde{M})) = \rk(H_0(M)) - 1$.
\par
Let $S' = S^x$. For $i \geq 0$, one has
\begin{equation*}
S'(i) = M(i) \cup \lbrace C_x(y) \colon y \in \tilde{M}(i-1) \rbrace.
\end{equation*}
Let $j\co C_i(M) \to C_i(S')$ be the map obtained from inclusion of $M(i)$ into $S'(i)$.
Let $h_i \co C_i(\tilde{M}) \to C_{i+1}(S')$ be the map defined by  $h_i([y]) = [C_x(y)]$.
One checks easily that
\begin{equation*}
 (h_{i-1} \partial + \partial h_i)([y] ) = j([y]).
\end{equation*}
(Sometimes we identify $C_{i}(M)$ as a subset of $C_i(S')$ via $j$ and write $[y]$ for $j([y])$).
Let 
\begin{equation*}
\bar{h}_i \co C_i(\tilde{M}) \to C_{i+1}(S')/ C_{i+1}(M)  
\end{equation*}
be the composition of $h_i$ with the projection map $C_{i+1}(S') \to C_{i+1}(S')/ C_{i+1}(M)$.
The map $\bar{h}_i$ is an isomorphism of abelian groups, since
$C_i(S') = C_i(M) \oplus h_{i-1}(C_{i-1}(\tilde{M}))$ for all $i \geq 0$. Moreover, the equation
$h \circ \partial + \partial \circ h = j$ shows that $(-1)^i \bar{h}_i$ is a chain isomorphism:
\begin{equation*}
 (-1)^i \bar{h}_i \co C_i(\tilde{M}) \simeq  C_{i+1}(S')/ C_{i+1}(M).
\end{equation*}
It follows that, there is an exact sequence of chain complexes,
\begin{equation*}
0 \to C_{\bullet}(M) \xrightarrow[]{j} C_{\bullet} (S') \xrightarrow[]{k} C_{\bullet - 1}(\tilde{M}) \to  0,
\end{equation*}
where, $k_i([z]) = 0$ for $z \in M_i$ and $k_i([C_x(y)]) = (-1)^i[y]$ for $y \in \tilde{M}_{i-1}$.
Taking the long exact sequence of homology groups, one gets
\begin{equation*}
\dotsb\to H_{i+1}(M) \xrightarrow[]{j_*} H_{i+1}(S') \xrightarrow[]{k_*} H_i(\tilde{M})
\xrightarrow[]{\delta_i} H_i(M) \to \dotsb
\end{equation*}
Let $\tau  = \sum_{\sigma} c_{\sigma} [\sigma] \in C_i(\tilde{M})$ be a $i$--cycle, that
is, $\partial \tau = 0$. The image of $\tau$ under the connecting
homomorphism $\delta_i$ is the homology class of $(j^{-1} \circ \partial \circ k^{-1}) (\tau)$,
where $k^{-1}(\tau)$ denotes any element in the pre-image. We have
\begin{equation*}
(-1)^i\tau = (-1)^i \sum c_{\sigma}[\sigma] \xrightarrow[]{k^{-1}}  \sum c_{\sigma} [C_x(\sigma)] \xrightarrow[]{\partial}
\tau' + \sum c_{\sigma} [\sigma] = \tau' + \tau,
\end{equation*}
where $\tau' = \sum_{z \in \tilde{M}_{i-1}} c'_{z} [C_x(z)]$ is a linear combination of ``cones''.
Since $k$ commutes with the boundary map, one has 
$k(\tau' + \tau) = k( \partial k^{-1} (-1)^i \tau) = \partial k k^{-1} (-1)^i \tau = \partial (-1)^i \tau =0$.
Since, by definition, $k$ ``kills" the old cells, one has $k(\tau) = 0$.
It follows that,
\begin{equation*}
0 = k (\tau' + \tau) = k(\tau') = \sum_{z \in \tilde{M}_{i-1}} (-1)^i c'_{z}  [z].
\end{equation*}
Thus $\tau' = 0$ and the connecting homomorphism $\delta_i\co H_i(\tilde{M}) \to H_i(M)$
is given by $\delta_i\co \tau \mapsto  (-1)^i\tau$.
 Since $H_i(M) = H_i(\tilde{M})$ for $i \geq 1$, the connecting homomorphism $\delta_i$
 is an isomorphism for $ i \geq 1$.
From the long exact sequence of homology groups, it follows that $H_i(S') = 0$ for $i > 1$.
It remains to calculate $H_1(S')$ and $H_0(S')$.
 \par
 For any $v \in M(0)$ one has $\partial[C_x(v)] = [v] - [C_x(\emptyset)]$,
implying that $[v]$ and $[C_x(\emptyset)]$ are in the same homology class in $H_0(S')$.
So $H_0(S') \simeq \ZZ$. Looking at the end of the long exact sequence, one has,
 \begin{equation*}
\dotsb \to H_1(\tilde{M}) \xrightarrow[]{\delta_1} H_1(M)
\to H_1(S') \to H_0(\tilde{M}) \to H_0(M) \to H_0(S') \to H_{-1}(\tilde{M}) \to 0.
 \end{equation*}
We know that $\delta_1$ is an isomorphism, $H_0(S') \simeq \ZZ$ and $H_{-1}(\tilde{M}) = 0$.
Using these informations, the above exact sequence reduces to
 \begin{equation*}
 0 \to H_1(S') \to H_0(\tilde{M}) \to H_0(M) \to \ZZ \to 0.
 \end{equation*}
But $H_0(\tilde{M})$ is a free $\ZZ$--module of rank one less than the rank of $H_0(M)$.
This forces $H_1(S') = 0$.
\end{proof}
\end{topic}
%
%

\begin{topic}
\label{baryviastellerpf}
\begin{proof}[proof of lemma \ref{baryviasteller}]
 The statements $\mathcal{A}(r+1)$, $\mathcal{B}(r)$ and $\mathcal{C}(r)$, for $0 \leq r \leq n$,
are proved by a single backward induction on $r$.
The last statement $(\mathcal{D})$ follows, by comparing the definition of the barycentric subdivision
$S^{(1)}$ with the description of $T_1$ provided by $\mathcal{A}(1)$ and $\mathcal{B}(0)$.
\par
To start induction, one has to check $\mathcal{A}(n+1), \mathcal{B}(n)$ and  $\mathcal{C}(n)$.
All these are obvious. The induction step goes as follows:
{\small
\begin{equation*}
\dotsb \implies \mathcal{A}(r+1) \implies \mathcal{B}(r) \implies \mathcal{C}(r) \implies
\mathcal{A}(r) \implies \mathcal{B}(r-1) \implies \mathcal{C}(r-1) \implies \dotsb
 \end{equation*}
}
Let $x \in T_r$. If $x \in T_{r+1}$ too, then we say that $x$ is an old cell of $T_r$. Otherwise,
we say that $x$ is a new cell of $T_r$.
\par
{\it proof of $\mathcal{B}(r)$ assuming $\mathcal{B}(m+1), \mathcal{C}(m+1), \mathcal{A}(m+1)$ for $m \geq r$: }
Suppose
\begin{equation*}
\beta = C_{w_j \dotsb w_1}(t) \leq C_{y_k \dotsb y_1} (v) = \alpha  \text{\;\; in \;\;}  T_{r+1}.
\end{equation*}
The cells of $T_{r+1}$ have this form because we are assuming $\mathcal{A}(r+1)$.
Next, $\mathcal{C}(r+1)$ implies that we can apply lemma \ref{commute} with
$T_{r+2} = X$ and $T_{r+1} = X_{(k)}$.
If both $\alpha$ and $\beta$ are old cells, then we are done by $\mathcal{B}(r+1)$.
If $\beta$ is old and $\alpha$ is new, then one must have
\begin{equation*}
C_{w_j \dotsb w_1}(t) \leq C_{y_{k-1} \dotsb y_1}(v) \text{\;\; in \;\;} T_{r+2} .
\end{equation*}
Now, $\mathcal{B}(r+1)$ implies that $\lbrace w_j < \dotsb < w_1 \rbrace$
is an ordered subset of $\lbrace y_{k-1} < \dotsb < y_1 \rbrace$
and $t \leq v$, from which we get $\mathcal{B}(r)$, in this case.
If $\beta$ is new, lemma \ref{commute} implies that $\alpha$ must also be new.
Further, one must have $y_k = w_j \in S(r+1)$ and
\begin{equation*}
C_{w_{j-1} \dotsb w_1}(t) \leq C_{y_{k-1} \dotsb y_1}(v) \text{\;\; in \;\;} T_{r+2}.
\end{equation*}
Using $\mathcal{B}(r+1)$, one gets, $\lbrace w_{j-1} < \dotsb < w_1 \rbrace$
is an ordered subset of $\lbrace y_{k-1} < \dotsb < y_1 \rbrace$ and $t \leq v$.
Together with $w_j = y_k$, the previous sentence implies $\mathcal{B}(r)$, in this case too.
\newline
\par
{\it proof of $\mathcal{C}(r)$ assuming $\mathcal{C}(m+1), \mathcal{A}(m+1), \mathcal{B}(m)$
for $m \geq r$: }
Suppose $ t \in S(r)$. From $\mathcal{A}(r+1)$, we know that $t \in T_{r+1}$. Suppose $t \leq C_{u_j \dotsb u_1}(v)$
in $T_{r+1}$. From $\mathcal{A}(r+1)$, it follows that $\rk(v) < r+1$, and from $\mathcal{B}(r)$, it follows that
$t \leq v$. But $\rk(t) = r$. So we must have $v = t$.
\newline
\par
{\it proof of $\mathcal{A}(r)$ assuming $\mathcal{A}(m+1), \mathcal{B}(m), \mathcal{C}(m)$ for $m \geq r$: }
Consider the transition from $T_{r+1}$ to $T_r$. The statement $\mathcal{A}(r+1)$ describes the
cells of $T_{r+1}$, while $\mathcal{B}(r)$ and $\mathcal{C}(r)$ describe the partial order on $T_{r+1}$. Let $t \in S(r)$. Note that $t$ ``survives'' as a cell of $T_{r+1}$.
One gets $T_r$ from $T_{r+1}$, by taking subdivision at each of these $t \in S(r)$. From
$\mathcal{C}(r)$, we know that the cells that ``die'' in the process of this subdivision are those
of the form $C_{u_j \dotsb u_1}(t)$, with $t \in S(r)$. So the old cells of $T_r$ are
\begin{equation*}
 \bigcup_{j = 0}^{n-r} \lbrace C_{u_j u_{j-1} \dotsb u_1}( v )  \colon
v < u_j < u_{j-1} \dotsb < u_1,  \rk(v) < r, \rk(u_{ j -i}) \geq r + 1 + i
\rbrace.
\end{equation*}
The new cells, that are ``born'' in this subdivision, have the form $C_t(x)$, where
$x \in M_{T_{r+1}}(t)$. Again, $\mathcal{C}(r)$ gives $x \leq C_{u_k \dotsb u_1}(t)$
in $T_{r+1}$. By $\mathcal{B}(r)$, this implies
$x = C_{w_j \dotsb w_1}(v)$, for some ordered subset,
$\lbrace w_j < \dotsb < w_1 \rbrace$, of $\lbrace u_k < \dotsb < u_1 \rbrace$
and some $v < t$. ($v = t$ is not a possibility, because $x \ngeq t$). In particular $w_{j -i} \geq u_{k -i}$.
It follows that there is a new  cell of $T_r$  of the form  $ C_t(C_{w_j \dotsb w_1}(v))$,
if and only if
\begin{equation*}
 v < t < w_j < \dotsb < w_1, \text{\;\;} \rk(t) = r, \text{\;\;} 0 \leq j \leq n - r, \text{\;\;}
\rk(w_{j-i}) \geq r + 1 + i,
\end{equation*}
where the last inequality follows from $ w_{j-i} \geq u_{k-i}$. The description of the cells of $T_r$ follows
by combining the descriptions of the old and the new cells.
\end{proof}
\end{topic}
%
%
\begin{topic}
\label{Phiphipf}
\begin{proof}[proof of equation \eqref{Phiphi} in \ref{baryinvf}]
We maintain the notations used in lemma \ref{baryviasteller}.
We can write $\varphi^{\circ}$ as a composition
$\varphi^{\circ} = \varphi^1 \circ \dotsb \circ \varphi^n$,
where $\varphi^r \co T_{r+1} \to T_r$ is the composite of the chain maps, 
given in \ref{StoSx}, corresponding to the stellar subdivisions
at the $r$--cells of $S$. Given $x_0 \in S(r)$ and a flag
$\gamma = \lbrace x_0 > x_1 > \dotsb > x_r \rbrace$
below $x_0$, we need to calculate the images of $[x_0]$ under successive application of
$\varphi^r$ and find the coefficient of
$[\gamma] = [C_{x_{r-1}  \dotsb x_1 x_0}(x_r)]$ in $\varphi^{\circ}([x_0])$.
\par
Clearly, $\varphi^{r+1} \circ \dotsb \circ \varphi^n [x_0] = [x_0] \in T_{r+1}$.
When we subdivide at the $r$--cells, $[x_0]$ is replaced by a linear combination
of cones (at the step where we take stellar refinement at $[x_0]$).
From the definition of the map $\varphi$, given in lemma \ref{StoSx}, we find that the image
of $[x_0]$ under $\varphi^r \circ \dotsb \circ \varphi^n$ is given by
\begin{equation*}
\varphi^r \circ \dotsb \circ \varphi^n [x_0] = \varphi^r([x_0]) = \sum_{x_1 \in \Delta x_0} s(x_0,x_1) [C_{x_0}(x_1)] \in T_r.
\end{equation*}
The statement $\mathcal{C}(r)$ in lemma \ref{baryviasteller} implies that
a cell of the form $C_{x_0}(x_1)$ dies at the next step, that is, during the transition
from $ T_r$ to $T_{r-1}$. More precisely, the cell $C_{x_0}(x_1)$
``dies", when we take stellar refinement at $x_1$. In that
step, $[C_{x_0}(x_1)]$ gets replaced by 
\begin{equation*}
\sum_{x_2 \in \Delta x_1} s(C_{x_0}(x_1),C_{x_0}(x_2)) [C_{x_1}(C_{x_0}(x_2))]
= -\sum_{x_2 \in \Delta x_1} s(x_1,x_2) [C_{x_1 x_0}(x_2)], 
\end{equation*}
(using equation \eqref{scone}). So 
\begin{equation*}
\varphi^{r-1} \circ \dotsb \circ \varphi^n [x_0] = -\sum_{x_1 , x_2} s(x_0,x_1) s(x_1,x_2)[C_{x_1 x_0}(x_2)] \in T_{r-1},
\end{equation*}
where the sum is over all $x_1$ and $x_2$ such that $x_1$ is a face of $x_0$ and $x_2$ is a face of $x_1$.
Continuing like this for $r$ steps, we find that
\begin{equation*}
\varphi^{\circ}[x_0] = 
\varphi^1 \circ \dotsb \circ \varphi^n [x_0] 
= \pm \sum_{x_1 \dotsb  x_r} \prod_{j=0}^{r-1}  s(x_j,x_{j+1}) [C_{x_{r-1} \dotsb x_0}(x_r)] \in T_1.
\end{equation*}
From \ref{defsigns} and our implicit assumption that $\omega_v(\lbrace v \rbrace) = 1$
for each zero--cell $v$, it follows that $\prod_{j=0}^{r-1}  s(x_j,x_{j+1}) = \omega_x(\gamma)$.
Thus, $\varphi^{\circ}$ matches the formula for $\Phi$ given in lemma \ref{StoS1}, up-to a sign.
\end{proof}
\end{topic}
%
%
%
%
%
%
%
%

%
%

\begin{thebibliography}{}
{\doublespacing
\bibitem[1]{PA:DR}
P.S. Alexandroff,
{\it Discrete Raume,}
Mathematiceskii Sbornik (NS) 2 (1937) 501--518.
%
\bibitem[2]{BD:HA}
J. Baez and J. Dolan,
{\it Higher-dimensional algebra III: $n$-categories and the algebra of opetopes,}
Adv. Math. 135(2) (1998) 145--206. 
%
\bibitem[3]{BM:SF}
J.A. Barmak and E.G. Minian, 
{\it Simple homotopy types and finite spaces,}
Adv. Math. 218 (2008) 87--104.
%
\bibitem[4]{CL:HC}
E. Cheng and A. Lauda,
{\it Higher-dimensional categories: an illustrated guide book,}
Preprint (2004). 
Available at http://www.dpmms.cam.ac.uk/\textasciitilde elgc2/guidebook 
%
\bibitem[5]{EZ:SS}
S. Eilenberg and J.A. Zilber,
{\it Semi--simplicial complexes and singular homology,}
Ann. of Math.(2) Vol. 51, No. 3, (1950) 499--513.

\bibitem[6]{ES:SS}
G. Ewald and G.C. Shephard,
{\it Stellar subdivisions of boundary complexes of convex polytopes,}
Math. Ann. 210 (1974) 7--16.
%
\bibitem[7]{GH:AG}
P. Griffiths and J. Harris,
{\it Principles of algebraic geometry,}
John Wiley \& Sons, Inc. Wiley Classics Library Ed. (1994).
%
\bibitem[8]{DK:CA}
D. Kozlov,
{\it Combinatorial Algebraic Topology,}
Springer Verlag, Series: Algorithms and computation in mathematics; 21 (2008).
%
\bibitem[9]{PM:SA}
J.P. May,
{\it Simplicial objects in algebraic topology,}
Chicago lectures in mathematics, Univ. of Chicago Press (1993).
%
\bibitem[10]{PM:FT}
J.P. May,
{\it Finite topological spaces,}
Notes for REU (2003). Available at
\newline
http://www.math.uchicago.edu/\textasciitilde may/MISCMaster.html
%
\bibitem[11]{PM:FS}
J.P. May,
{\it Finite spaces and simplicial complexes,}
Notes for REU (2003) Available at
\newline
http://www.math.uchicago.edu/\textasciitilde may/MISCMaster.html 
%
\bibitem[12]{MM:SH}
M.C. McCord,
{\it Singular homology groups and homotopy groups of finite
topological spaces,}
Duke Mathematical Journal 33 (1966) 465--474.
%
\bibitem[13]{RS:AS}
R. Street,
{\it The algebra of oriented simplexes,}
J. Pure Appl. Algebra 49(3) (1987) 283--335.
%
\bibitem[14]{MW:PT}
M. Wachs,
{\it Poset Topology: Tools and Applications,}
Geometric combinatorics, IAS/Park City Math. Ser. 13, Amer. Math. Soc. Providence, RI, (2007) 497--615.
%
\bibitem[15]{JW:SH}
J.H.C. Whitehead,
{\it Simple homotopy types,}
Amer. J. Math. 72 (1950) 1--57.

}
%
\end{thebibliography}
\end{document}